\documentstyle[12pt,amssymb,amscd]{amsart}  

\newtheorem{thm}{Theorem}[subsection]

\newtheorem{cor}[thm]{Corollary}
\newtheorem{prop}[thm]{Proposition}
\newtheorem{lemma}[thm]{Lemma}

\newtheorem{definition}[thm]{Definition}
\newtheorem{proposition}[thm]{Proposition}

\theoremstyle{remark}
\newtheorem{remark}[thm]{Remark}
\newtheorem{example}[thm]{Example}

\theoremstyle{definition}
\newtheorem{defn}[thm]{Definition}

\numberwithin{equation}{section}

\newcommand{\C}{{\Bbb C}}
\newcommand{\bbA}{{\Bbb A}}
\newcommand{\bbC}{{\Bbb C}}

\newcommand{\bbL}{{\Bbb L}}

\newcommand{\bbW}{{\Bbb W}}
\newcommand{\bbZ}{{\Bbb Z}}

\newcommand{\cO}{{\cal O}}

\newcommand{\cM}{{\cal M}}

\newcommand{\cD}{{\cal D}}

\newcommand{\cA}{{\cal A}}

\newcommand{\cI}{{\cal I}}

\newcommand{\cE}{{\cal E}}

\newcommand{\cT}{{\cal T}}

\newcommand{\One}{{\boldmath 1}}

\newcommand{\R}{\operatorname{\bold R}}
\newcommand{\Hom}{\operatorname{Hom}}
\newcommand{\isomo}{\overset{\sim}{=}}
\newcommand{\isomoto}{\overset{\sim}{\to}}
\newcommand{\id}{\operatorname{id}}

\newcommand{\shHom}{\underline{\operatorname{Hom}}}

\newcommand{\Symm}{\operatorname{Sym}}
\newcommand{\g}{{\frak{g}}}

\thanks
{The first and the third authors were supported
in part by NSF grants.}

\begin{document}

\title{Riemann-Roch theorems via deformation quantization II}

\author{P.Bressler}
\address{Institut des Hautes Recherches Scientifiques, Bur-sur-Yvette, France}
\email{bressler@@ihes.fr}
\author{R.Nest}
\address{Mathematical Institute, Universitatsparken, 5,
2100 Copenhagen, Denmark}
\email{rnest@@math.ku.dk}
\author{B.Tsygan}
\address{Department of Mathematics, The Pennsylvania State University,
University Park, PA 16802, USA}
\email{tsygan@@math.psu.edu}

\maketitle

\section{Introduction}
In this paper we give a proof of a Riemann-Roch type theorem for symplectic deformations.
The main result of the present paper is used in \cite{BNT1} in order to establish a
conjecture of P.~Schapira and J.-P.~Schneiders about the microlocal Euler class.

 Let ${\cE _X}$ be the sheaf of micro-differential operators on a complex manifold $X$.
The microlocal Euler class $\mu \operatorname{eu}(\cM ^{\bullet})$ of a good complex
$\cM ^{\bullet}$ of $\cD$-modules on $X$ is the image of a certain homology class
$\operatorname{Eu} ^{\cE _X} (\cM ^{\bullet})$  under the trace density morphism
$\mu ^{\cE _X}$ from the Hochschild complex of ${\cE _X} $ to the de Rham complex of
$T^{\ast}X$. We showed in \cite{BNT1} that $\operatorname{Eu} ^{\cE _X} (\cM ^{\bullet})$
extends to a homology class $\operatorname {ch}^{\cE _X} (\cM ^{\bullet})$ of the periodic
cyclic complex, and that the map $\mu ^{\cE _X}$ extends to a morphism $\tilde{\mu} ^{\cE _X}$
from the periodic cyclic complex of $\cE _X$ to the de Rham complex of $T^*X$. Next we replaced the algebra of microdifferential
operators by a deformation quantization $\bbA _{T^{\ast}X}^t$ of the structure sheaf of the
cotangent bundle of $X$. We defined in this context the trace density map $\tilde {\mu} ^{\bbA ^t}$
and the periodic cyclic homology class $\operatorname {ch}^{\bbA ^t} (\cM ^{\bullet})$.
This allowed us to reduce the conjecture of Schapira and Schneiders to a theorem about deformation
quantizations of the structure sheaf $\cO _M$ of a complex manifold $M$ with a holomorphic symplectic
form $\omega$. Namely, let $\bbA ^t _M$ be such a deformation. We define the trace density map 
\begin{equation} \label{eq:muccper}
\tilde{\mu}^t : CC^{\operatorname{per}}_{\bullet}(\bbA^t _M) \rightarrow  \bbC _M [[t,t^{-1}][[u,u^{-1}] 
\end{equation}
in the derived category of sheaves on $M$. Here $CC^{\operatorname{per}}_{\bullet}(\bbA^t _M)$ is the
periodic cyclic complex of the sheaf of algebras  $\bbA^t _M$ and $u$ is a formal parameter of degree $2$.
We also defined the principal symbol map 
\begin{equation} \label{eq:muccpero}
\tilde{\mu} \circ \sigma : CC^{\operatorname{per}}_{\bullet}(\bbA^t _M) @>>>
CC^{\operatorname{per}}_{\bullet}(\cO _M) @>>> \bbC _M [[u,u^{-1}]
\end{equation}
(the first map is reduction modulo $t$, the second is the Connes-Hochschild-Kostant-Rosenberg
map from \cite{C}, \cite{L}). We reduced the conjecture of Schapira and Schneiders to a Riemann-Roch
theorem for deformation quantizations: 
\begin{equation} \label{eq:rrtat}
\tilde{\mu}^t = (\tilde {\mu} \circ \sigma) \smile \widehat{A}(TM) \smile e^{\theta}
\end{equation}
\begin{equation} \label{eq:muccpero1}
H^{\bullet}(M,CC^{\operatorname{per}}_{\bullet}(\bbA^t _M)) @>>>
H^{\bullet}(M, \bbC _M [[u,u^{-1}][[t,t^{-1}])
\end{equation}

(theorem \cite{BNT1}, Theorem 4.6.1).

Here $\theta \in \frac{1}{t}H^2(M, \bbC)[[t]]$ is the characteristic class of the deformation
(\cite{D}, \cite{F}; cf. the end of Section \ref{defqua}). 

Note that the trace density map generalizes the canonical trace from \cite{F}; it is also related to the construction of \cite{CFS} (cf. \cite{Ha} for explanation).

 To prove such a theorem, we need two techniques:

\begin{itemize}
\item{deformation quantization of complex manifolds}
\item{cyclic complexes of algebras}
\end{itemize}
as well as 
\begin{itemize}
\item{Lie algebra cohomology and Gelfand-Fuks map}
\end{itemize}

In section \ref{defqua} we review deformation quantization of complex symplectic manifolds (i.e.
symplectic manifolds with a holomorphic symplectic structure). We essentially follow \cite{NT3}
where results of Fedosov are adapted to the complex case more systematically than here. We define
Fedosov connections and their gauge equivalences. We show that locally any two are gauge equivalent.
This implies that the sheaf of horizontal sections of a Fedosov connection is a deformation quantization
of the structure sheaf of algebras. Finally we show that any deformation of the structure sheaf comes
from a Fedosov connection.

In section 4 we define in our context the Gelfand-Fuks map which allows to reduce problems about any
deformation quantization to problems about Lie algebra cohomology. Along these lines we reduce Riemann-Roch
theorem to a parallel statement about complexes of Lie algebra cochains (theorem \ref{thm:formalRR}). This
theorem compares two morphisms 
\begin{equation} \label{twomor}
C^{\bullet}({\frak{g}}, {\frak{h}}; \;CC^{\text{per}}_{\bullet}(W)) @>>>
C^{\bullet}({\frak{g}}, {\frak{h}}; {\widehat{\Omega}}^{\bullet}({\Bbb{C}}^{2d}))[[t,t^{-1}][[u, u^{-1}]
\end{equation}
where $W$ is the Weyl algebra of the complex symplectic space ${\Bbb{C}}^{2d}$, ${\frak{g}}$ is the Lie
algebra of continuous derivations of $W$, ${\frak{h}}$ is the subalgebra ${{\frak{sp}}}(2d)$ of
${\frak{g}}$, ${\widehat{\Omega}}^{\bullet}$ stands for the completion of the de Rham complex of
$\C ^{2d}$ at the origin, $C^{\bullet}({\frak{g}}, {\frak{h}}; \;\ldots)$ is the relative Lie algebra cochain
complex with coefficients in a complex of $({\frak{g}}, {\frak{h}})$-modules, and $CC^{\text{per}}_{\bullet}$ 
is the periodic cyclic complex. The two morphisms (\ref{twomor}) are the trace density map and the principal
symbol map which are defined in \ref{chclie}. A little more precisely, the trace density morphism maps the periodic
cyclic complex of $W$ to the complex 
\begin{equation} \label{eq:mirrorsymp}
({\widehat{\Omega}}^{\bullet}[2d][[t,t^{-1}][[u, u^{-1}], td_{DR})
\end{equation}
and the principal symbol morphism maps the periodic cyclic complex to
\begin{equation} \label{eq:mirrorclass}
({\widehat{\Omega}}^{-\bullet}[[t,t^{-1}][[u, u^{-1}], ud_{DR})
\end{equation}
We identify both complexes above with the complex
\begin{equation} \label{eq:mirror}
({\widehat{\Omega}}^{\bullet}[2d][[t,t^{-1}][[u, u^{-1}], d_{DR})
\end{equation}
The complex (\ref{eq:mirrorsymp}) is identified with (\ref{eq:mirror}) by the operator $I$ which is equal to the multiplication operator by $t^{i-d}$ on ${\widehat{\Omega}}^{i}[[t,t^{-1}][[u, u^{-1}]$; the complex (\ref{eq:mirrorclass}) is identified with (\ref{eq:mirror}) by the operator $J$ which is equal to the multiplication operator by $u^{i-d}$ on ${\widehat{\Omega}}^{i}[[t,t^{-1}][[u, u^{-1}]$. Note the symmetry between "the Planck constant" $t$ and the Bott periodicity formal parameter $u$.
Thus, Riemann-Roch theorem is reduced to a purely algebraic statement (theorem \ref{thm:formalRR}).

To prove this theorem, we construct an operator on $CC^{\text{per}}_{\bullet}(W))[[t, t ^{-1}]$ whose composition with the trace density map is easily computable and which is cohomologous to the identity map. To be able to do that, we use an algebra of operations on the periodic cyclic complex of any algebra \cite{NT2}. More precisely, we construct an action of the reduced cyclic complex $\overline{C}^{\lambda}_{\bullet - 1}(A)$ on $CC^{\text{per}}_{\bullet}(A)$ for any algebra $A$ over a commutative ring $k$ of characteristic zero (Theorem \ref{thm:calc3}). At the level of homology this action is very simple: it is given by multiplication by the image of the boundary map
\begin{equation} \label{boumap}
\partial: \overline{C}^{\lambda}_{\bullet}(A) \to CC_{\bullet}(k)[1]
\end{equation}
(the right hand side is equal to $k[u^{-1},u]]/k[[u]]$).
At the level of chains, however, the above action is given by a non-trivial formula which allows to compute it explicitly in an important example. Namely, for $A=W$ we construct the fundamental class in $\overline{C}^{\lambda}_{2d - 1}(W)$ which we extend, canonically up to  homotopy, to a cocycle $U$ of the double complex $C^{\bullet}({\frak{g}}, {\frak{h}}; \overline{C}^{\lambda}_{ \bullet}(W))$. One sees easily that the pairing with $U$, followed by the trace density map, is the principal symbol map. It remains to compute the pairing with $U$ at the level of homology. This is , as mentioned above, tantamount to computing the image of $U$ in $H^{\bullet}({\frak{g}}, {\frak{h}}; CC_{\bullet}(k))  \simeq H^{\bullet}({\frak{g}}, {\frak{h}}) \otimes (k[u^{-1},u]]/k[[u]])$
  under the boundary map $\partial$.

For this, two more technical tools are needed. First, we need an explicit formula for the boundary map $\partial$ (\ref{boumap})). We recall this formula, due to Brodzki, in subsection \ref{ss:computu}. This formula, being completely explicit, does not, however, solve our problem by itself because we lack explicit information about the fundamental cycle $U$. To be able to write $U$ explicitly, we need one more simplification. Recall first that for a Lie algebra ${\frak{g}}$ with a subalgebra  ${\frak{h}}$ whose adjoint action on ${\frak{g}}$ is semisimple, one can define the Chern-Weil map
$$ {\operatorname{CW}}: k[{\frak{h}}^{\ast}]^{\frak{h}} \to C^{\bullet}({\frak{g}}, {\frak{h}}) $$
Composed with the Gelfand-Fuks map (\ref{map:GF}), this map gives the classical Chern-Weil map. Following \cite{NT4}, we generalize the map ${\operatorname{CW}}$ to the case of cochains with non-trivial coefficients. More precisely, we consider complexes $L^{\bullet}$ of ${\frak{g}}$-modules which are homotopy constant in a simple natural way which is explained in Definition \ref{hotri}. For such complexes we construct a generalized Chern-Weil map whose right hand side is the double complex $C^{\bullet}({\frak{g}}, {\frak{h}};L^{\bullet})$. The left hand side is, roughly speaking, the complex of $L^{\bullet}$-valued ${\frak{h}}$-equivariant polynomials.  If ${\frak{g}}$ acts on a manifold $M$ by vector fields, then the de Rham complex $L^{\bullet}=\cA ^{\bullet}(M)$ is a homotopy constant complex of ${\frak{g}}$-modules. If the action of ${\frak{h}}$ integrates to a group action, then the right hand side of the generalized Chern-Weil map is the complex of equivariant different
ial forms on $M$. 

(Note that the formula of Brodzki for the boundary map itself is in our exposition obtained from the map dual to  $\operatorname{CW}$, applied in the case when ${\frak{g}}={\frak{gl}}(A)$ and ${\frak{h}}={\frak{gl}}(k)$ for any algebra $A$. The Chern-Weil map surfaces in this paper one more time, when we deduce theorems \ref{thm:calc2} and \ref{thm:calc3}
from a more general theorem of non-commutative differential calculus). 

We show that the fundamental cocycle $U$ is in the image under the generalized Chern-Weil map of a cocycle that can be constructed explicitly. The image of this last cocycle under the boundary map $\partial$ is computed directly. This completes the proof.

Let us say a few words about one of the tools used in the proof, namely operations on Hochschild and cyclic complexes (subsection \ref{ss:opera}). The results and constructions of these subsections are part of what we call non-commutative differential calculus: non-commutative generalizations of geometric objects such as differential forms, multivector fieds, etc., carry algebraic structures similar to classical ones, up to strong homotopy. Another manifestation of this principle is a recent theorem of Tamarkin \cite{T} stating that the algebra of Hochschild cochains is a strong homotopy Gerstenhaber algebra. Conjecturally, the theorem of Tamarkin and the results of 
\ref{ss:opera} have far-reaching common generalizations. Together with the formality theorem of Kontsevich, they would lead to a generalization of the Riemann-Roch theorem to all deformations, symplectic or not (\cite{TT}, \cite{Ts1}).

{\bf Acknowledgements}. This paper, as well as \cite{BNT1}, has profited greatly from our discussions with A.~D'Agnolo, M.~Kashiwara and P.~Schapira. We, especially the last author, are greatly indebted to B.~Feigin for many ideas used in this work. 

\section{Deformation quantization}
\label{section:DQ}
In this section we outline the approach to deformation quantization pioneered
by Fedosov which is crucial to our method of proof of Theorem 4.6.1, \cite{BNT1}.

\subsection{Review of deformation quantization}  \label{defqua}
A {\em deformation quantization} of a manifold $M$ is a formal one parameter
deformation of the structure sheaf $\cO_M$, i.e. a sheaf of algebras
$\bbA^t_M$ flat over $\bbC [[t]]$ together with an isomorphism of
algebras $\psi :\bbA^t_M\otimes_{\bbC [[t]]}\bbC {\rightarrow}\cO_M$.

The formula
\[
\lbrace f, g\rbrace = \frac{1}{t}[\tilde f,\tilde g ]
+ t\cdot \bbA^t_M\ ,
\]
where $f$ and $g$ are two local sections of $\cO_M$ and $\tilde f$, $\tilde g$
are their respective lifts to $\bbA^t_M$, defines a Poisson structure on $M$
called the Poisson structure associated to the deformation quantization
$\bbA^t_M$.

A deformation quantization $\bbA^t_M$ is called {\em symplectic} if the 
associated Poisson structure is nondegenerate. In this case $M$ is symplectic.
In what follows we will only consider symplectic deformation quantizations,
so assume that $\bbA^t_M$ is symplectic from now on.

Let us note first that, given a deformation ${\mathbb A}^t_M$ as above, $\psi$
induces locally an isomorphism of the sheaves of ${\Bbb{C}}$-vector spaces:
$$
\tilde{\psi }: \prod (t^{n} {\mathbb A}^t_U )
/t^{n+1} {\mathbb A}^t_U \simeq {\cal O}_U [[ t ]].
$$
If one choses a cover by Stein subsets, then the left hand side is isomorphic to ${\mathbb A}^t_U$; this implies that there exist local
isomorphisms
$$
\Phi_U :{\mathbb A}^t_U \rightarrow  ({\cal O}_U [[ t ]] ,*_{U} )
$$
of ${\Bbb{C}}[[t]]$-algebras with transition isomorphisms $G_{UV} = \Phi_U \Phi_V^{-1}$ of
the form
\begin{equation}
G_{UV} = f +t D_1^{UV} +t^2 D_2^{UV} +\ldots .  
\end{equation}

In the rest of this section we will work under following condition.

{\bf {{\underline {Assumption.}} Both the local products $*_U$ and the linear transformations $D_i^{UV}$ above
are given by finite order holomorphic (bi-)differential operators}}. 

Similarly, an isomorphism of two deformations will be always understood as an isomorphism of sheaves of algebras whose restrictions to coordinate charts $U$ are given by 
$$T_U (f) = \operatorname{id} + \sum_{i \geq 1} t^i T_{U, i}(f)$$
where $T_{U, i}(f)$ are holomorphic differential operators. 

\begin{example} \label{weylal}
Let $M={\bbC^{2d}}$.
Define the product on $\bbA^t_{\bbC^{2d}}$ to be, in coordinates
$x_1,\ldots,x_d,\xi_1,\ldots,\xi_d$ on $\bbC^{2d}$, the standard
Moyal--Weyl product
\[
(f\ast g)(\underline x,\underline\xi) = \\
exp\left(\frac{t}{2}\sum_{i=1}^d\left(
\frac{\partial\ }{\partial\xi_i}\frac{\partial\ }{\partial y_i}-
\frac{\partial~}{\partial\eta_i}\frac{\partial~}{\partial x_i}\right)\right)
f(\underline x,\underline\xi)g(\underline y,\underline\eta)
\vert_{\overset{\underline x=\underline y}{\underline\xi =\underline\eta}}
\]
where $\underline x = (x_1,\ldots,x_d),\ \underline\xi = (\xi_1,\ldots,\xi_d),
\ \underline y = (y_1,\ldots,y_d),\ \underline\eta = (\eta_1,\ldots,\eta_d)$.
\end{example}

\subsection{The Weyl algebra}
The Weyl algebra of a symplectic vector space $(V,\omega)$ over $\bbC$
may be considered as a symplectic deformation quantization of
the completion of $V$ at the origin.

Here we briefly recall the definition and the basic properties of the
Weyl algebra $W = W(V,\omega)$ of a (finite dimensional) symplectic vector space 
$(V,\omega)$ over $\bbC$. Let $V^\ast = \Hom_{\bbC}(V,\bbC)$.

The Moyal-Weyl product on $\Symm^\bullet(V^\ast)[t]$ is defined
by the formula
\begin{equation}\label{formula:MW}
f\ast g = \sum_{n=0}^\infty\frac{1}{n!}\left(\frac{t}{2}\right)^n
\omega(d^nf,d^ng)\ ,
\end{equation}
where
\[
d^n : \Symm^\bullet(V^\ast) @>>>
\Symm^\bullet(V^\ast)\otimes \Symm^n(V)
\]
assigns to a polynomial the symmetric tensor composed of its
$n$-th order partial derivatives, and $\omega$ is extended naturally to
a bilinear form on $\Symm^n(V)$. Clearly, this sum is finite for any two polynomials $f$, $g$.  Of course, (\ref{formula:MW})  becomes the formula from \ref{weylal} in Darboux coordinates.

Let 

$${\frak{m}} = {\operatorname{Ker}}(\Symm^\bullet(V^\ast) \rightarrow {\Bbb{C}})$$

\[
\widehat\Symm^\bullet(V^\ast) = \varprojlim
\frac{\Symm^\bullet (V^\ast)}{{\frak{m}}^n}
\]
Define a filtration $F_{\bullet}$ on $\widehat\Symm^\bullet(V^\ast)[[t]]$ by
\begin{equation} \label{filtra}
F_{-p}\widehat\Symm^\bullet(V^\ast)[[t]] = \sum_{i+\frac{j}{2} \geq p}{\widehat{\frak{m}}}^i {t}^j
\end{equation}

The Moyal-Weyl product
endows $\widehat\Symm^\bullet(V^\ast)[[t]]$ with a structure of an
associative algebra with unit over $\bbC[[t]]$. Moreover, the Moyal-Weyl product is
continuous in the topology defined by the filtration $F_{\bullet}$. Let $W=W(V)$ denote
the topological algebra over $\bbC[[t]]$ whose underlying
$\bbC[[t]]$-module is $\widehat\Symm^\bullet(V^\ast)[[t]]$ and the
multiplication is given by the Moyal-Weyl product.

Note also that the center of $W$ is equal to $\bbC[[t]]$ and
$\left[W,W\right] = t\cdot W$.

Clearly, the association $(V,\omega)\mapsto W(V)$ is functorial. In
particular, the group $Sp(V)$ acts naturally on $W(V)$ by continuous
algebra automorphisms. The infinitesimal action of the Lie algebra ${\frak{sp}}(V)$ is given by operators $\frac{1}{t} \operatorname{ad}(f)$ where $f$ are in $\Symm ^2 (V^{\ast})$. This correspondence defines a linear isomorphism ${\frak{sp}}(V) \simeq \Symm ^2 (V^{\ast})$.

\subsection{Derivations of the Weyl algebra}
Let ${\frak{g}} = {\frak{g}}(V)$ denote the Lie algebra of continuous,
$\bbC[[t]]$-linear derivations of $W$. Then there is a central extension
of Lie algebras
\begin{equation} \label{eq:ext}
0 @>>> \frac{1}{t}\bbC[[t]] @>>>
\frac{1}{t}W @>>> {\frak{g}} @>>> 0\ ,
\end{equation}
where the Lie algebra structure on $\frac{1}{t}W$ is
given by the commutator (note that
$\left[\frac{1}{t}W,\frac{1}{t}W\right]
\subseteq\frac{1}{t}W$) and the second map is defined
by $\frac{1}{t}f\mapsto \frac{1}{t}\left[f,\bullet
\right]$.

Let
\[
F_p{\frak{g}} = \left\{D\in{\frak{g}}\ \vert\ D(F_iW)\subseteq F_{i+p}W
\ \text{for all $i$}\right\}
\]
Then $({\frak{g}},F_\bullet)$ is a filtered Lie algebra and the action
of ${\frak{g}}$ on $W$ respects the filtrations, i.e.
$\left[F_p{\frak{g}},F_q{\frak{g}}\right]\subseteq F_{p+q}{\frak{g}}$ and
$F_p{\frak{g}}F_qW\subseteq F_{p+q}W$.

The following properties of the filtered Lie algebra $({\frak{g}},F_\bullet)$
are easily verified:
\begin{enumerate}
\item $Gr^F_p{\frak{g}} = 0$ for $p>1$ (in particular ${\frak{g}}=F_1{\frak{g}}$),
hence $Gr^F_1{\frak{g}}$ is Abelian;

\item the composition $\frac{1}{t}V\hookrightarrow
\frac{1}{t}W\to{\frak{g}}\to Gr^F_1{\frak{g}}$ is an
isomorphism;

\item the composition ${{\frak{sp}}}(V)\to\Symm^2(V^\ast)\to
\frac{1}{t}F_{-2}W \to F_0{\frak{g}}\to Gr^F_0{\frak{g}}$
is an isomorphism;

\item under the above isomorphisms the action of $Gr^F_0{\frak{g}}$ on
$Gr^F_1{\frak{g}}$ is identified with the natural action of
${\frak{sp}}(V)$ on $V$ (in particular there is an isomorphism
${\frak{g}}/F_{-1}{\frak{g}}\isomo V\ltimes{\frak{sp}}(V)$);

\item the Lie algebra $F_{-1}{\frak{g}}$ is pro-nilpotent.
\end{enumerate}

In what follows ${\frak{h}}$ will denote the Lie subalgebra of ${\frak{g}}$
which is the image of the embedding ${\frak{sp}}(V)\hookrightarrow{\frak{g}}$.

\subsection{The Weyl algebra as a deformation quantization}
The (commutative) algebra $W/t W$ is naturally isomorphic to
the completion $\widehat\cO =\widehat\cO_V$ of the ring $\cO_V$ of regular functions
on $V$ at the origin, i.e. with respect to the powers of the maximal ideal
$\frak{m}$ of functions which vanish at $0\in V$. The natural surjective
map
\[
\sigma : W(V) @>>> \widehat\cO_V
\]
is strictly compatible with the filtration $F_{\bullet}$ on $W$ and the
$\frak{m}$-adic filtration on $\widehat\cO_V$. The Poisson bracket induced on $\widehat {\cO} _V$ is the standard one. 

It is well known (and not difficult to show) that any deformation of $\widehat\cO_V$ with the above properties is isomorphic to $W$. An isomorphism can be chosen as a continuous linear map which preserves the filtration $F$ and is equal to identity on $\operatorname {gr}_F {\mathbb A}^t$. 

The Lie algebra ${\frak{g}}$ acts by derivations on $\widehat\cO$
by the formula
\[
D(f) = \sigma(D(\tilde f))\ ,
\]
where $D\in {\frak{g}}$ and $\tilde f\in W$ is such that $\sigma(\tilde f)=f$.
Thus, $\sigma$ is a map of ${\frak{g}}$-modules.
\subsection{The sheaf of Weyl algebras}
The sheaf of Weyl algebras $\bbW_M$ on $M$ is the sheaf of topological
algebras over the sheaf of topological algebras $\cO_M[[t]]$ (equipped
with the $t$-adic topology) defined as follows.

Let $\Theta_M$ denote the sheaf of holomorphic vector fields on $M$, $\Theta_M ^{\ast} = Hom_{\cO_M} (\Theta_M, {\cO_M})$.
Let $\cI$ denote the kernel of the augmentation map
\[
\Symm_{\cO_M}^\bullet(\Theta_M ^{\ast})\otimes_{\cO_M}\cO_M
@>>> \cO_M\ .
\]
Let $\widehat\Symm_{\cO_M}^\bullet(\Theta_M ^{\ast})$ be the completion $\widehat\Symm_{\cO_M}^\bullet(\Theta_M ^{\ast})[[t]]$
of $\Symm_{\cO_M}^\bullet(\Theta_M ^{\ast})\otimes_{\cO_M}\cO_M$
in the $\cI$-adic topology. Then $\widehat\Symm_{\cO_M}^\bullet(\Theta_M ^{\ast})[[t]]$ is a topological $\cO_M[[t]]$-module. The
Weyl multiplication on $\widehat\Symm_{\cO_M}^\bullet(\Theta_M ^{\ast})[[t]]$
is defined by \eqref{formula:MW}.
Then $\bbW_M$ is the sheaf of $\cO_M[[t]]$-algebras whose underlying
sheaf of $\cO_M[[t]]$-modules is 
$\widehat\Symm_{\cO_M}^\bullet(\Theta_M ^{\ast})[[t]]$ and the multiplication
is given by the Moyal-Weyl product. Note that the center of $\bbW_M$ is 
$\cO_M[[t]]$, and $\left[\bbW_M,\bbW_M\right] = t\cdot\bbW_M$.

Let 
\begin{equation}
F_{-p}\bbW_M = \sum_{i+\frac{j}{2} \geq p} \widehat\cI^{i} t ^j
\end{equation}
 Then $(\bbW_M,F_\bullet)$ is
a filtered ring, i.e. $F_p\bbW_M\cdot F_q\bbW_M\subseteq F_{p+q}\bbW_M$.
Note that the quotients $F_p\bbW_M/F_q\bbW_M$ are locally free
$\cO_M$-modules of finite rank. The Weyl multiplication is continuous in the topology induced by the filtration $F$.
\subsection{Review of the Fedosov construction}
We refer the reader to \cite{F} and \cite{NT3} for a detailed exposition of
the construction of deformation quantizations via Fedosov connections.

Let $\cA^{p,q}_M$ denote the sheaf of complex valued $C^\infty$-forms of
type $(p,q)$ on $M$, $\cA^s_M = \oplus_{p+q=s}\cA^{p,q}_M$. Let
$F_\bullet\cA^\bullet_M$ denote the Hodge filtration and $d$ the de Rham
differential. Then $((\cA^\bullet_M,d),F_\bullet)$ is a filtered differential
graded algebra (i.e. $F_r\cA^\bullet_MF_s\cA^\bullet_M\subseteq
F_{r+s}\cA^\bullet_M$ and $d(F_r\cA^\bullet_M)\subseteq F_r\cA^\bullet_M$).

\begin{remark} \label{rmk:tens}
Here and below, for any ring $A$, when we consider a tensor product over $A$ by a pro-finite rank free $A$-module, we mean a completed tensor product. 
\end{remark}
The space
\[
\cA^\bullet_M(\bbW_M)\overset{def}{=}\bigoplus_p \cA^p_M(\bbW_M)[-p]
\]
has a natural structure of a sheaf of graded algebras. Let
\begin{equation} \label{filfor}
F_r\cA^\bullet_M(\bbW_M) = \sum_{p+q=r}F_p\cA^\bullet_M
\otimes_{\cO_M}F_p\bbW_M\ .
\end{equation}
Then $(\cA^\bullet_M(\bbW_M),F_\bullet)$ is a filtered graded algebra over
$(\cA^\bullet_M,F_\bullet)$ (i.e.\linebreak
$F_r\cA^p_M(\bbW_M)F_s\cA^q_M(\bbW_M)\subseteq
F_{r+s}\cA^{p+q}_M(\bbW_M)$ and $F_r\cA^p_MF_s\cA^q_M(\bbW_M)\subseteq
F_{r+s}\cA^{p+q}_M(\bbW_M)$).

The associated graded algebra of $\cA^\bullet_M(\bbW_M)$ of the filtration $F_\bullet$ is the commutative algebra
\begin{equation} \label{grafor}
\cA^\bullet_M (\Symm_{\cO_M}^\bullet(\Theta_M ^{\ast}))[t] = {\cA}^{0,\bullet}_M (\Symm_{\cO_M}^\bullet(\Theta_M ^{\ast}) \otimes \wedge_{\cO_M}^\bullet(\Theta_M ^{\ast}))[t]
\end{equation}
(the degree of an element of ${\cA}^{0,\bullet}_M (\Symm^i V\otimes \wedge ^j V)t ^k$ is $-(i+j+2k)$).
\begin{defn}  \label{delta}
Let $\delta$ be the Koszul differential on (\ref{grafor}), i.e. the ${\cA}^{0,\bullet}_M (\wedge_{\cO_M}^\bullet(\Theta_M ^{\ast}))[t]$-linear derivation of the commutative algebra which sends any $X \in \Symm ^1 (\Theta_M ^{\ast})$ to $X \in \wedge ^1 (\Theta_M ^{\ast})$
\end{defn}
We will show that for any deformation of $\cO _M$ there exists a map
\begin{equation}  \label{nabla}
\nabla : \cA^\bullet_M(\bbW_M) @>>> \cA^\bullet_M(\bbW_M)[1]
\end{equation}
which has the following properties:
\begin{enumerate}
\item $\nabla\left(F_p\cA^\bullet_M(\bbW_M)\right)\subseteq
F_p\cA^\bullet_M(\bbW_M)$; the induced map
\[
Gr^F_{\bullet}\nabla : Gr^F_{\bullet}\cA^\bullet_M(\bbW_M) @>>> 
Gr^F_{\bullet}\cA^\bullet_M(\bbW_M)[1]
\]
is equal to $\overline{\partial} + \delta$ where $ \delta$ is the Koszul differential from Definition \ref{delta};
\item $\nabla^2 = 0$;

\item $(\cA^\bullet_M(\bbW_M),\nabla)$ is a sheaf of differential graded 
algebras over the $C^\infty$ de Rham complex $(\cA^\bullet_M,d)$
(in particular $H^\bullet(\cA^\bullet_M(\bbW_M),\nabla)$ is a sheaf of
graded algebras over the constant sheaf $\bbC_M[[t]]$);

\item there is a filtered quasi-isomorphism
\[
(\bbA^t_M,F_\bullet) @>>> ((\cA^\bullet_M(\bbW_M),\nabla),F_\bullet)
\]
of differential graded algebras over 
 $\bbC_M[[t]]$ (in particular  \newline
$H^p(F_q\cA^\bullet_M(\bbW_M),\nabla)=0$ for $p\neq 0$)
\end{enumerate} 
Therefore $\nabla$ is determined by its component
\[
\nabla : \cA^0_M(\bbW_M) @>>> \cA^1_M(\bbW_M)
\]
which has all the properties of a flat connection on $\bbW_M$.

To be able to prove the existence and the properties of $\nabla$, let us define it in a slightly different, more differential-geometric language. 

Let $H = Sp(\dim M)$ and ${\frak{h}}$ the Lie algebra of $H$. Put $P @>{\pi}>> M$ denote the $H$-principal bundle
of symplectic frames in $TM$. We identify $TM$ with the vector bundle associated
to the standard representation of $H$. For $h$ in ${\frak{h}}$ let $\tilde h$ be the vector field on $P$ induced by $h$.

Recall that, for an ${\frak{h}}$-module $V$, one can define the subcomplex
$\left[\pi_*\cA_P^\bullet\otimes V\right]^{basic}\subset\pi_*\cA_P^\bullet$
of all basic $V$-valued forms $\alpha$, i.e. of all forms in $\cA_P^\bullet\otimes V$ such that 
$$\iota_{\tilde{h}} \alpha = 0; \; L_{\tilde{h}} \alpha = h\alpha$$
for all $h$ in $\frak{h}$. In other words, the complex of basic forms can be defined as the pull-back diagram
\[
\begin{CD}
\left[\pi_*\cA_P^\bullet\otimes V\right]^{basic} @>>>
C^\bullet(\pi_*\cT_P,{\frak{h}};V) \\
@VVV  @VVV \\
\pi_*\cA_P^\bullet\otimes V @>>> C^\bullet(\pi_*\cT_P;V)
\end{CD}
\]
where $\cT_P$ denotes the sheaf of Lie algebras of $C^\infty$ vector fields
on $P$ and the complexes in the right hand side stand for the relative and the absolute complexes of Lie algebra cochains.

Let $W$ denote the Weyl algebra
of the standard representation of $H$. Then $\bbW_M$ is identified with
the sheaf of sections of the associated bundle $P\times_HW$ and
pull-back by $\pi$
\[
\cA_M^\bullet(\bbW_M) @>>> \pi_*\cA_P^\bullet \otimes W
\]
identifies $\bbW_M$-valued forms on $M$ with the subcomplex of basic
$W$-valued forms on $P$. The flat connection $\nabla$ gives rise to
a ${\frak{g}}(=Der(W))$-valued 1-form
$A\in H^0(P;\cA_P^1 \otimes{\frak{g}})$ which satisfies 
\begin{equation}\label{cond:conn-form}
\text{$\iota_{\tilde h}A=h$ and $L_{\tilde h}A = ad(h)(A)$ for all
$h\in{\frak{h}}$}
\end{equation}
(where, as above, $\tilde h$ is the vector field on $P$ induced by $h$)
and the Maurer-Cartan equation 
\begin{equation}  \label{eq:MC}
dA+\frac12\lbrack A,A\rbrack = 0
\end{equation}
 (so that
$(d + A)^2=0$). Then, pull back by $\pi$ induces the isomorphism of
filtered complexes
\[
(\cA_M^\bullet(\bbW_M),\nabla) @>>>
\left(\left\lbrack\pi_*\cA_P^\bullet\widehat\otimes W\right\rbrack^{basic},
d+A\right)\ .
\]
In other words, let $U_\alpha$ be a cover of M by Darboux coordinate charts. The coordinates on $U_{\alpha}$ trivialize the holomorphic tangent bundle on those charts, and let $g_{\alpha \beta}$ be the transition isomorphisms with respect to these trivializations. We have holomorphic maps
$$ g_{\alpha \beta}: U_{\alpha} \cap U_{\beta} \rightarrow Sp(2n,{\Bbb{C}})  \rightarrow \operatorname{Aut} (W)$$
which are the transition functions of the Weyl bundle. A $\bbW_M$-valued form on $M$ is a collection of $W$-valued forms $\omega _{\alpha}$ satisfying 
$\omega _{\alpha} = g_{\alpha \beta}\omega _{\beta}$. A form $A$ as above is a collection of ${\frak{g}}$-valued 1-forms $A_{\alpha}$ satisfying
$$A_{\alpha} = -dg_{\alpha \beta}\;g_{\alpha \beta}^{-1} + g_{\alpha \beta} A_{\beta} g_{\alpha \beta}^{-1}$$

Any ${\frak{g}}$-valued 1-form $A$ on $P$ which satisfies
(\ref{eq:MC}) defines a map $\nabla$ as in (\ref{nabla}). 
This map automatically satisfies conditions 2 and 3 after (\ref{nabla}). In order for it to satisfy condition 1, the following must be true:
\begin{equation} \label{eq:fedos1}
A \in A_{-1} + A_0 + \cA ^1 (M, F_{-1}{\frak{g}})
\end{equation}
\begin{equation} \label{eq:fedos2}
A_0 \in \cA ^{1,0} ({P}, {\frak{sp}}(2n))
\end{equation}
and $A_{-1}$ is the tautological holomorphic form 
\begin{equation} \label{eq:fedos3}
A_{-1}: T M \isomoto T^{\ast} M \to F_1 \operatorname{Der}(W(T))
\end{equation}
where the first isomorphism is induced by $- \omega$ viewed as an element of $ T^*M {\otimes ^2}\simeq \operatorname{Hom}(TM, T^*M)$.

\begin{defn}\label{dfn:Fedosov}
A ${\frak {g}}$-valued form $A$ on $P$ which satisfies (\ref{cond:conn-form}), (\ref{eq:MC}), (\ref{eq:fedos1}), (\ref{eq:fedos2}), (\ref{eq:fedos3}) is called a Fedosov connection form. The operator $\nabla = d+A$ is called a Fedosov connection.
\end{defn}
\begin{defn}\label{defini:Fedosov}
Two Fedosov connections $d+A$ and $d+B$ are gauge equivalent if there is an element $X$ in $\cA ^0 (M, F_{-1} {\frak {g}})$ such that $d+B = \operatorname{exp} (\operatorname{ad}X)(d+A)$
\end{defn}

\begin{thm} \label{thm:fed}
\begin{enumerate}

\item For any Fedosov connection $\nabla$, the sheaf $\operatorname{Ker}\nabla$ is isomorphic to a deformation of $\cO_M$.
\item For any deformation of $\cO_M$ there is a Fedosov connection $\nabla$ and a quasi-isomorphism $\bbA^t_M @>>> (\cA^\bullet_M(\bbW_M),\nabla)$
of (sheaves of) differential graded algebras over $\bbC_M[[t]]$. If two deformations $\operatorname{Ker}\nabla$ and $\operatorname{Ker}\nabla '$ are isomorphic then $\nabla$ and $\nabla ' $ are gauge equivalent. 
\item Any two Fedosov connections are locally gauge equivalent.
\end{enumerate}
\end{thm}

For the readers convenience we include a sketch of the proof. A more complete proof is contained in \cite{NT3}.

{\underline{\bf {Sketch of the proof. }}}
To see why any two Fedosov connections are locally gauge equivalent, assume that we are given two Fedosov connections $ \nabla$ and $\nabla '$. We proceed by induction, assuming that $\nabla - \nabla '$ is in $\cA^{\bullet}(M, F_{-n}{\frak{g}})$, $n\geq 0$. Then there is a gauge transformation which sends $\nabla '$ to $\nabla ^{''}$ such that $\nabla - \nabla ^{''}$ is in $\cA^{1}(M, F_{-n-1}{\frak{g}})$. Indeed, to find such a transformation of the form $\operatorname{exp} (\operatorname{ad}X)$ where $X$ is in $\cA^{0}(M, F_{-n-1}{\frak{g}})$, one has to solve the equation
\begin{equation} \label{equququ}
 (\overline{\partial} + \delta) (X) = \nabla - \nabla ^{''} modulo \cA^{1}(M, F_{-n-1}{\frak{g}})
\end{equation}
(cf. condition 1 after (\ref{nabla})). 

Consider the complex $(\cA^{\bullet}_M ( gr^{F}_{\bullet}({\frak{g}})), \overline{\partial} + \delta )$. It is isomorphic to $(\cA^{\bullet}_M( \wedge ^{\bullet}(T^*_M) \otimes \widehat{S} ^{\bullet}(T^*_M)[[t]]), \overline{\partial} + \delta _0 )$ where $\delta _0$ is induced by the de Rham differential  on $\wedge ^{\bullet}(T^*_M) \otimes \widehat{S} ^{\bullet}(T^*_M)[[t]]$. Therefore $\delta$ is acyclic in degree greater than zero, and its cohomology is isomorphic to the Dolbeault complex $\cA^{\bullet}_M[[t]]$. Thus, the complex  $(\cA^{\bullet}_M ( gr^{F}_{\bullet}({\frak{g}})), \overline{\partial} + \delta )$ is quasi-isomorphic to 
$\cA^{\bullet}_M[[t]]$.

One checks that the right hand side of (\ref{equququ}) is a $\overline{\partial} + \delta$-cocycle in  $\cA^{\bullet}_M ( gr^{F}_{-n}({\frak{g}})$. Therefore, locally, this equation always has a solution.

 To prove that any deformation is isomorphic to ${\operatorname{Ker}}\nabla$ for some Fedosov connection $\nabla$, note that, if $\bbA^t _M$ is such a deformation, then the bundle of jets of its holomorphic sections is a bundle of algebras. This buldle of algebras is equipped with a canonical flat connection which we call the Grothendieck connection. It also has a filtration which is preserved by this connection. To define this filtration on any fiber $J_x$ of the jet bundle at a point $x \in M$, let ${\frak{m}}_x$ be the ideal of jets of functions vanishing at $x$, and let $F^n J_x$ be the completion of 
$$F^n J_x = \sum_{i+\frac{j}{2} \geq n} t^j {\frak{m}}_x^i$$
 The associated graded algebra of this filtration is the Weyl algebra. As $C^{\infty}$ bundles of algebras, the bundle of jets and its (completed) associated graded bundle are isomorphic. Therefore the Weyl bundle has a flat connection (the image of the canonical connection under this isomorphism). This connection is a Fedosov connection whose kernel is isomorphic to the sheaf of holomorphic sections of $\nabla$.  

To compute the cohomology of the complex $(\cA _M ^{\bullet}(\bbW), \nabla)$, note that, locally, any Fedosov connection is gauge equivalent to the standard connection $\nabla _0 = \overline{\partial} + \delta$ whose cohomology is the sheaf $\cO [[t]]$. So, for any $\nabla$, the zeroth cohomology is a deformation of $\cO_M$ and the higher cohomology is zero. Finally, any isomorphism ${\operatorname{Ker}}\nabla \simeq {\operatorname{Ker}}\nabla '$ extends uniquely to an automorphism of the Weyl bundle intertwining $\nabla$ and $\nabla '$. (To prove this, it is enough to work locally and to assume that $\nabla = \nabla ' = \nabla _0$). This completes the proof.

\subsection{The characteristic class of a symplectic deformation} \label{ss:theta}
To a deformation $\bbA^t _M$, one can associate a characteristic class $\theta $ in 
\newline
$H^2 (M, \frac{1}{t}\bbC[[t]])$ as follows. Let $A$ be a Fedosov connection such that $\bbA^t _M = {\operatorname{Ker}}(d+A)$. Let $\tilde{A}$ be a $\frac{1}{t}W$-valued form on $P$ satisfying (\ref{cond:conn-form}), (\ref{eq:fedos1}), (\ref{eq:fedos2}), (\ref{eq:fedos3}). We call $\tilde{A}$ a lifting of $A$ if its image in $\cA^{1}(M, {\frak{g}})$ under the map induced by (\ref{eq:ext}).

It is easy to see that any Fedosov connection has a lifting. For any lifting, the curvature 
\begin{equation} \label{eq:theta}
\theta = d\tilde{A} + \frac{1}{2}[\tilde{A}, \tilde{A}]
\end{equation} 
is a 2-form with values in $\frac{1}{t}\bbC[[t]]$ because of (\ref{eq:MC}). If $\tilde{A}$ and $\tilde{A} + \alpha$ are two liftings of the same connection, then $\alpha$ is $\bbC[[t]]$-valued and the curvatures of the two liftings differ by $d\alpha$. For two gauge equivalent liftings the curvatures are the same. Therefore, because of theorem \ref{thm:fed} (second part of statement 2), the cohomology class of $\theta$ depends only on an isomorphism class of a deformation. This cohomology class is called the characteristic class of the deformation $\bbA^t _M$.

\section{The Riemann-Roch formula in the formal setting}
\label{section:formalRR}

The aim of this section is to state Theorem \ref{thm:formalRR} which compares two morphisms from the periodic cyclic complex of the Weyl algebra $W$ to the de Rham complex.  The first morphism is the trace density map (cf. \ref{tracedensity}) which maps the periodic cyclic complex of $W$ to the complex \begin{equation} \label{eq:mirrorsymp1}
({\widehat{\Omega}}^{\bullet}[2d][[t,t^{-1}][[u, u^{-1}], td_{DR})
\end{equation}
The second is the the principal symbol map followed by the Hochschild-Kostant-Rosenberg map (cf. \ref{sssection:HKR})  which maps the periodic cyclic complex to
\begin{equation} \label{eq:mirrorclass1}
({\widehat{\Omega}}^{-\bullet}[[t,t^{-1}][[u, u^{-1}], ud_{DR})
\end{equation}
We identify both complexes above with the complex
\begin{equation} \label{eq:mirror1}
({\widehat{\Omega}}^{\bullet}[2d][[t,t^{-1}][[u, u^{-1}], d_{DR})
\end{equation}
The complex (\ref{eq:mirrorsymp1}) is identified with (\ref{eq:mirror1}) by the operator $I$ which is equal to $t^{i-d}$ on ${\widehat{\Omega}}^{i}[[t,t^{-1}][[u, u^{-1}]$; the complex (\ref{eq:mirrorclass1}) is identified with (\ref{eq:mirror1}) by the operator $J$ which is equal to $u^{i-d}$ on ${\widehat{\Omega}}^{i}[[t,t^{-1}][[u, u^{-1}]$. Riemann-Roch theorem for deformation quantizations will be deduced from Theorem \ref{thm:formalRR} in Section
\ref{section:GF}.

\subsection{The Hochschild homology of the Weyl algebra}
\label{subsection:HH-W}
We recall the calculation of the Hochschild homology of the Weyl algebra
(\cite{FT1}, \cite{Bry}).

From now on $W=W(V,\omega)$. Note that the identity map furnishes the
identification $W^{op}\isomo W(V,-\omega)$. Let
\[
W^{\frak{e}}\overset{def}{=}W((V,\omega)\oplus (V,-\omega))\ .
\]
There is an isomorphism of algebras $W\widehat\otimes_{\bbC[[t]]}W^{op}
\rightarrow W^{\frak{e}}$.

Let
\begin{equation}  \label{eq:C(W)}
{C}_p(W)\overset{def}{=}W((V,\omega)^{\oplus p+1})\ .
\end{equation}
Here, as above, for us tensor products of complete $k$-modules always mean completed tensor products.

The Hochschild differential \eqref{diffl:b} extends to this setting and
gives rise to the complex ${C}_\bullet(W)$ which represents
$W\otimes^{\bbL}_{W^{\frak{e}}}W$ in the derived category.

The Koszul complex
$(K^\bullet,\partial)$ is defined by
\[
K^{-q} = W\otimes{\bigwedge}^q V^\ast
\]
with the differential acting by
\begin{equation*}
\partial(f\otimes v_1\wedge\ldots\wedge v_q) = t
\sum_i(-1)^i [f,v_i]\otimes v_1\wedge\ldots\wedge\widehat v_i\wedge
\ldots\wedge v_q
\end{equation*}
Here we consider $V^\ast$ embedded in $W$ and $\bigwedge^qV^\ast$ embedded
in $(V^\ast)^{\otimes q}$.

The map
\[
K^\bullet @>>> C_\bullet(W)
\]
defined by
\[
f\otimes v_1\wedge\ldots\wedge\ldots\wedge v_q
\mapsto f\otimes Alt(v_1\otimes\cdots\otimes v_q)
\]
is easily seen to be a quasi-isomorphism. (Indeed, it preserves the complete filtration by powers of $t$ and induces a quasi-isomorphism at the level of associated graded quotients).
Hence it induces a quasi-isomorphism
\begin{equation}\label{map:K-to-C}
K^\bullet[t^{-1}]
@>>> C_\bullet(W)[t^{-1}]\ .
\end{equation}

Suppose that $\dim V = 2d$.
Let $\widehat\Omega^\bullet =\widehat\Omega^\bullet_V$ denote the de Rham
complex of $V$ with
formal coefficients (i.e. $\widehat\Omega^q_V =\Omega^q_V\otimes_{\cO_V}
\widehat\cO_V$). 
The map
\[
W\otimes{\bigwedge}^q V^\ast @>>>
\widehat\Omega^{2d-q}[[t]]
\]
given by
\[
f\otimes v_1\wedge\ldots\wedge v_q \mapsto
f\cdot\iota_{v_1}\cdots\iota_{v_q}(\omega^{\wedge d})
\]
(where $f\in\widehat\cO$)
is easily seen to determine an isomorphism of complexes
\begin{equation}\label{map:K-to-DR}
(K^\bullet(W),\partial)
@>>> (\widehat\Omega^\bullet[[t]], t\cdot d_{DR})[2d]\ .
\end{equation}

It follows from Poincar\'{e} Lemma (for power series) that the following is true:
\begin{proposition} \label{prop:HH(W)}
\begin{enumerate}

\item $H^p{C}_\bullet(W)[t^{-1}] = 0$ for $p\neq -2d$ 
\item the maps \eqref{map:K-to-C} and \eqref{map:K-to-DR} induce the
isomorphism
$H^{-2d}{C}_\bullet(W)[t^{-1}]\isomo\bbC[t^{-1},t]]$
which maps the class of the cycle
$1\otimes\frac{1}{d!}\omega^{\wedge d}\in W\otimes
{\bigwedge}^{2d}V^\ast$ to $1$.
\end{enumerate}
\end{proposition}
If $x_1,\ldots ,x_d,\xi_1,\ldots ,\xi_d$ is (dual to) a symplectic basis
of $V$ (so that $\omega = \sum_i x_i\wedge\xi_i$), then  one can choose a generator which 
is represented by the cycle $Alt(1\otimes x_1\otimes\cdots\otimes x_d\otimes
\xi_1\otimes\cdots\otimes\xi_d)\in C_{2d}(W)$. We will denote this cycle
(and its class) by $\Phi = \Phi_V$. Note also that $\Phi$ corresponds
under the above (quasi)isomorphisms to the cocycle $1\in\widehat\Omega^0$.

Observe that the Lie algebra ${\frak{g}}$ acts on the Hochschild complex by Lie derivatives. The cycle $\Phi$ is not invariant under the action
of ${\frak{g}}$. It is, however, invariant under the action of the subalgebra
${\frak{h}}$.


\subsection{Some characteristic classes in Lie algebra cohomology} \label{chclie}
We will presently construct the classes in relative Lie algebra cohomology
of the pair $({\frak{g}},{\frak{h}})$ which enter the Riemann-Roch formula in
the present setting.

\subsubsection{The trace density}  \label{tracedensity}
In \ref{subsection:HH-W}, we constructed a quasi-isomorphism
\[
(\widehat\Omega^\bullet[2d][t^{-1},t]], td_{DR}) @>>> C_\bullet(W)[t^{-1}]
\]
Since ${\frak{h}}$ acts semi-simply on $C_\bullet(W)[t^{-1}]$ and
$\widehat\Omega^\bullet[t^{-1},t]]$, this quasi-isomorphism admits an ${\frak{h}}$-equivariant splitting
\[
\mu^t_{(0)} :
C_\bullet(W)[t^{-1}] @>>> (\widehat\Omega^\bullet_V[2d][t^{-1},t]], td_{DR})
\]
which is a quasi-isomorphism. We identify the right hand side with $(\widehat\Omega^\bullet_V[2d][t^{-1},t]], d_{DR})$ by the operator $I$ as in the beginning of Section \ref{section:formalRR}.

We will consider the map $\mu^t_{(0)}$ as a relative Lie algebra cochain
\[
\mu^t_{(0)}\in
C^0({\frak{g}},{\frak{h}};\Hom^0(C_\bullet(W)[t^{-1}],
\widehat\Omega^\bullet[2d][t^{-1},t]]))\ .
\]

\begin{lemma}
$\mu^t_{(0)}$ extends to a cocycle $\mu^t = \sum_p\mu^t_{(p)}$
with
\[
\mu^t_{(p)}\in
C^p({\frak{g}},{\frak{h}};\Hom^{-p}(C_\bullet(W)[t^{-1}],
\widehat\Omega^\bullet[t^{-1},t]][2d]))\ .
\]
Moreover, any two such extensions are cohomologous.
\end{lemma}

The proof follows easily from proposition \ref{prop:HH(W)}.

\begin{remark} \label{rmk:hom}
As above, when we work with a ${\frak{g}}$-module $C_p (W)$, $C_p(\widehat{\cO }_V)$, etc., we always use completed tensor products. By ${\operatorname{Hom} }(C^p(W), \widehat{\Omega}^q [t^{-1},t]])$, etc. we always mean the space of those maps from $\otimes ^{p+1} \widehat{\cO}_V [[t]]$ to $\widehat {\Omega}^q [t^{-1},t]]$ which are of the form 
$$D(f_0, \ldots, f_p) = \sum _{i \geq m} t^i D_i ( f_0, \ldots, f_p)$$
where $D_i$ are multi-differential operators (of finite order) which are equal to zero if $f_j = 1$ for some $j > 0$.
\end{remark}

Cup product with $\mu^t$ induces the quasi-isomorphism of complexes
\[
\mu^t : C^\bullet({\frak{g}},{\frak{h}};C_\bullet(W)[t^{-1}]) @>>>
C^\bullet({\frak{g}},{\frak{h}};\widehat\Omega^\bullet[2d][t^{-1},t]])
\]
unique up to homotopy.

\begin{lemma}
$\mu^t_{(0)}$ extends to an ${\frak{h}}$-equivariant quasi-isomorphism of
complexes
\[
\tilde\mu^t_{(0)} : CC^{per}_\bullet(W)[t^{-1}] @>>>
\widehat\Omega^\bullet[2d][t^{-1},t]][u^{-1},u]]
\]
\end{lemma}

We will consider the map $\tilde\mu^t_{(0)}$ as a relative Lie algebra
cochain
\[
\tilde\mu^t_{(0)}\in
C^0({\frak{g}},{\frak{h}};\Hom^0(CC^{per}_\bullet(W)[t^{-1}],
\widehat\Omega^\bullet[2d][t^{-1},t]][u^{-1},u]]))\ .
\]

\begin{lemma}\label{lemma:formalTR}
$\tilde\mu^t_{(0)}$ extends to a cocycle $\tilde\mu^t =
\sum_p\tilde\mu^t_{(p)}$
with
\[
\mu^t_{(p)}\in
C^p({\frak{g}},{\frak{h}};\Hom^{-p}(CC^{per}_\bullet(W)[t^{-1}],
\widehat\Omega^\bullet[2d][t^{-1},t]][u^{-1},u]]))\ .
\]
Moreover, any two such extensions are cohomologous.
\end{lemma}

The proof follows easily from proposition \ref{prop:HH(W)}.

\begin{remark}  \label{rmk:cchom}
In the formula above, $\Hom$ is understood as the module of those homomorphisms which preserve the $(u)$-adic filtration and are $\bbC[[t,t^{-1}][[u]]$-linear.
\end{remark} 

\begin{remark} \label{rmk:derived}
The image of $\mu ^t$ under the canonical morphism from the complex \newline
$C^\bullet({\frak{g}},{\frak{h}};\Hom^\bullet
(CC^{per}_\bullet(W)[t^{-1}],\widehat\Omega^\bullet[2d][t^{-1},t]]))$ to \newline
$\R\Hom^\bullet_{({\frak{g}},{\frak{h}})}(CC^{per}_\bullet(W)[t^{-1}],
\widehat\Omega^\bullet[2d][t^{-1},t]][u^{-1},u]])$
 a well defined morphism
\[
\tilde\mu^t :
CC^{per}_\bullet(W)[t^{-1}] @>>>
\widehat\Omega^\bullet[2d][t^{-1},t]][u^{-1},u]]
\]
in the derived category of $({\frak{g}},{\frak{h}})$-modules. The image of
$\tilde\mu^t$ under the functor of forgetting the module structure
is $\tilde\mu^t_{(0)}$.
\end{remark}

Cup product with $\tilde\mu^t$ induces the quasi-isomorphism of complexes
\[
\tilde\mu^t : C^\bullet({\frak{g}},{\frak{h}};CC^{per}_\bullet(W)[t^{-1}])
@>>>
C^\bullet({\frak{g}},{\frak{h}};
\widehat\Omega^\bullet[2d][t^{-1},t]][u^{-1},u]])
\]
unique up to homotopy.

The natural inclusion
\[
\iota : CC^{per}_\bullet(W) @>>> CC^{per}_\bullet(W)[t^{-1}]
\]
is a morphism of complexes of $({\frak{g}},{\frak{h}})$-modules, therefore
determines a cocycle
\[
\iota\in C^0({\frak{g}},{\frak{h}};\Hom^0
(CC^{per}_\bullet(W), CC^{per}_\bullet(W)[t^{-1}]))\ .
\]

The cup product of $\iota$ and $\tilde\mu^t$ is a cocycle
\[
\tilde\mu^t\smile\iota\in
C^\bullet({\frak{g}},{\frak{h}};\Hom^\bullet
(CC^{per}_\bullet(W),\widehat\Omega^\bullet[2d][t^{-1},t]][u^{-1},u]])
\]
of (total) degree zero which represents the morphism
\[
\tilde\mu^t\circ\iota :
CC^{per}_\bullet(W) @>>> \widehat\Omega^\bullet[2d][t^{-1},t]][u^{-1},u]]
\]
in the derived category of $({\frak{g}},{\frak{h}})$-modules.
\subsubsection{The symbol and the Hochschild-Kostant-Rosenberg map}
\label{sssection:HKR}
The symbol map is the map
\[
\sigma : CC^{per}_\bullet(W) @>>> CC^{per}_\bullet(\widehat\cO)
\]
induced by the algebra homomorphism $\sigma : W\to\widehat\cO$). The Hochschild-Kostant-Rosenberg map
\[
\tilde\mu : CC^{per}_\bullet(\widehat\cO) @>>> 
(\widehat\Omega^\bullet[u^{-1},u]] , ud_{DR}) @>>>
(\widehat\Omega^\bullet[u^{-1},u]][t^{-1},t]],  ud_{DR})
\]
is defined by $f_0\otimes\cdots\otimes f_p\mapsto\frac{1}{p!}
f_0df_1\wedge\ldots\wedge df_p$. We identify the last complex with 
\begin{equation} \label{eq:mirror2}
(\widehat\Omega^\bullet[2d][u^{-1},u]][t^{-1},t]],  d_{DR})
\end{equation}
by the operator $J$ as in the beginning of Section \ref{section:formalRR}. 

Both $\tilde{\mu}$ and $\sigma$ are morphisms of complexes of $({\frak{g}},{\frak{h}})$-modules, therefore they
determine cocycles
\[
\sigma\in C^\bullet({\frak{g}},{\frak{h}};\Hom^\bullet
(CC^{per}_\bullet(W), CC^{per}_\bullet(\widehat\cO)))
\]
and
\[
\tilde\mu\in C^\bullet({\frak{g}},{\frak{h}};\Hom^\bullet
(CC^{per}_\bullet(\widehat\cO),
\widehat\Omega^\bullet [2d][t^{-1},t]][u^{-1},u]]))\ .
\]

\subsubsection{The $\widehat A$-class}\label{sssection:formalA}
Recall from \ref{chclasstrivco} the Chern-Weil map
\[
CW : \widehat\Symm^\bullet({\frak{h}}^\ast)^{{\frak{h}}} @>>>
C^{2\bullet}({\frak{g}},{\frak{h}};\bbC)
\]

Let $\widehat A$ denote the image under the Chern-Weil map of
\[
{\frak{h}}\ni X\mapsto\det\left(\frac{ad(\frac{X}{2})}{exp(ad(\frac{X}{2})) -
exp(ad(-\frac{X}{2}))}\right)\ .
\]

\subsubsection{The characteristic class of the deformation}
\label{sssection:charcl}
The central extension of Lie algebras
\[
0 @>>> \frac{1}{t}\bbC[[t]] @>>>
\frac{1}{t}W @>>> {\frak{g}} @>>> 0
\]
restricts to a trivial extension of ${\frak{h}}$, therefore is classified by
a class $\theta\in H^2({\frak{g}},{\frak{h}};\frac{1}{t}\bbC[[t]])$
represented by the cocycle
\[
\theta : X\wedge Y\mapsto \widetilde{[X,Y]} -
[\widetilde{X},\widetilde{Y}]
\]
where $\widetilde{(\ )}$ is a choice of a $\bbC[[t]]$-linear splitting
of the extension.

It will be important to remember that the above cohomology class is also a partial case of the Chern-Weil map:
\[
CW : \frac{1}{t}\widehat\Symm^\bullet(\tilde{\frak{h}}^\ast)^{\tilde{\frak{h}}} @>>>
C^{2\bullet}(\tilde{\frak{g}},\tilde{\frak{h}};\frac{1}{t}\bbC [[t]]) @>>> C^{2\bullet}({\frak{g}},{\frak{h}};\frac{1}{t}\bbC[[t]])
\]

\subsection{Riemann-Roch formula in Lie algebra cohomology}
Let 
$A = \sum_{k\geq 0} A_{2k}$ be an even cohomology class, $A_{2k} \in H^{2k}({\frak{g}},{\frak{h}})$. Let $L^{\bullet}$ is any complex of ${\frak{g}}[u^{-1},u]]$-modules.By $A\smile ?$ we denote the multiplication operator on $C^\bullet({\frak{g}},{\frak{h}}; L^{\bullet})$ given by 
\begin{equation} \label{eq:smile}
c \mapsto A\smile c = \sum_{k\geq 0} A_{2k}\cdot c\cdot u^{-k}
\end{equation}
where $c$ is a cochain $C^\bullet({\frak{g}},{\frak{h}}; L^{\bullet})$. Formula (\ref{eq:smile}) makes $A\smile ?$ an operator of degree zero on $C^\bullet({\frak{g}},{\frak{h}}; L^{\bullet})$ (recall that degree of $u$ is equal to $2$).

The following theorem is the analog of Theorem 4.6.1 of \cite{BNT1} in the
setting of this section.

\begin{thm}\label{thm:formalRR}
The cocycle $\tilde\mu^t\circ\iota - (\widehat A \cdot
e^\theta)\smile(\tilde\mu \circ \sigma) $ is cohomologous to zero in
$C^\bullet({\frak{g}},{\frak{h}};\Hom^\bullet
(CC^{per}_\bullet(W),
\widehat\Omega^\bullet [2d][t^{-1},t]][u^{-1},u]]))$.
\end{thm}

\begin{cor}
The diagram
\[
\begin{CD}
CC^{per}_\bullet(W) @>{\sigma}>> CC^{per}_\bullet(\widehat\cO) \\
@V{\iota}VV @VV{(\widehat A\cdot e^\theta)\smile \tilde\mu}V \\
CC^{per}_\bullet(W)[t^{-1}] @>{\tilde\mu^t}>>
\widehat\Omega^\bullet_V [2d][t^{-1},t]][u^{-1},u]]
\end{CD}
\]
in the derived category of $({\frak{g}},{\frak{h}})$-modules is commutative.
\end{cor}

\begin{cor}
The diagram
\[
\begin{CD}
C^\bullet({\frak{g}},{\frak{h}};CC^{per}_\bullet(W)) @>{\sigma}>>
C^\bullet({\frak{g}},{\frak{h}};CC^{per}_\bullet(\widehat\cO)) \\
@V{\iota}VV
@VV{(\widehat A\cdot e^\theta)\smile \tilde\mu}V
\\
C^\bullet({\frak{g}},{\frak{h}};CC^{per}_\bullet(W)[t^{-1}])
@>{\tilde\mu^t}>> C^\bullet({\frak{g}},{\frak{h}};
\widehat\Omega^\bullet [2d][t^{-1},t]][u^{-1},u]])
\end{CD}
\]
is homotopy commutative.
\end{cor}

A proof of Theorem \ref{thm:formalRR} will be given in Section \ref{section:formalRR-pf}.

\section{Gelfand-Fuks map}\label{section:GF}
In this section we introduce the machinery related to the
Gelfand-Fuks map and reduce Theorem 4.6.1, \cite{BNT1} from its formal analog
Theorem \ref{thm:formalRR}.

Suppose given a complex manifold $M$ and a symplectic deformation
quantization $\bbA^t_M$ of $M$. By theorem \ref{thm:fed}, there exist, unique up to gauge equivalence, a Fedosov connection $\nabla$ such that $\bbA^t_M \isomoto \operatorname{Ker}\nabla$. Let $\omega\in H^0(M;\Omega^2_M)$ 
denote the associated symplectic form.

Given $A$ as above and a (filtered) topological ${\frak{g}}$-module $L$ such
that the action of ${\frak{h}}\subset{\frak{g}}$ integrates to an action of $H$.
Set
\[
(\cA_M^\bullet(L),\nabla)\overset{def}{=}
\left(\left\lbrack\pi_*\cA_P^\bullet\widehat\otimes L\right\rbrack^{basic},
d+A\right)\ .
\]

Note that the association $L\mapsto (\cA_M^\bullet(L),\nabla)$ is functorial
in $L$. In particular it extends to complexes of ${\frak{g}}$-modules.

Taking $L=\bbC$, the trivial ${\frak{g}}$-module, we recover
$(\cA_M^\bullet,d)$. For any complex $(L^\bullet,d_L)$ of ${\frak{g}}$-modules
as above the complex $(\cA_M^\bullet(L^\bullet),\nabla + d_L)$ has a natural
structure of a differential graded module over $(\cA_M^\bullet,d)$.

The Gelfand-Fuks map
\[
GF = GF_{\nabla} : C^\bullet({\frak{g}},{\frak{h}};L) @>>> \cA_M^\bullet(L)
\]
(the source understood to be the constant sheaf) is defined by the
formula
\[
GF(c)(X_1,\ldots,X_p) = c(A(X_1),\ldots,A(X_p))\ ,
\]
where $c\in C^p({\frak{g}},{\frak{h}};L)$ and $X_1,\ldots,X_p$ are locally
defined vector fields. It is easy to verify that $GF$ takes values in
basic forms and is a map of complexes. Note also that $GF$ is natural
in $L$. In particular the definition above has an obvious extension to
complexes of ${\frak{g}}$-modules. Let $g = \operatorname{exp}X$ be a gauge transformation between $\nabla$ and $\nabla '$, Then $g$ induces an isomorphism 
$$g_{\ast} : (\cA_M^\bullet(L),\nabla) \isomoto \cA_M^\bullet(L),\nabla' )$$
It is easy to see that $GF_{\nabla '}$ and $g_{\ast} \circ GF_{\nabla}$ are canonically homotopic. 

We now proceed to apply the above constructions to particular examples
of complexes $L^\bullet$ of ${\frak{g}}$-modules. In all examples below
the ${\frak{g}}$-modules which appear have the following additional
property which is easy to verify, namely,
\[
\text{$H^p(\cA_M^\bullet(L),\nabla) = 0$ for $p\neq 0$}\ .
\]
If $(L^\bullet,d_L)$ is a complex of ${\frak{g}}$-modules with the above
property, then the inclusion
\begin{equation}\label{map:incl-ker}
(\ker(\nabla),d_L\vert_{\ker(\nabla)})\hookrightarrow
(\cA_M^\bullet(L^\bullet),\nabla + d_L)
\end{equation}
is a quasi-isomorphism. In such a case the map \eqref{map:incl-ker} induces
an isomorphism
\[
\R\Gamma(M;\ker(\nabla)) @>>> \Gamma(M;\cA_M^\bullet(L^\bullet))
\]
in the derived category of complexes since the sheaves
$\cA_M^\bullet(L^\bullet)$ are fine and the Gelfand-Fuks map induces
the (natural in $L^\bullet$) morphism
\begin{equation}\label{map:GF}
GF : C^\bullet({\frak{g}},{\frak{h}};L) @>>> \R\Gamma(M;\ker(\nabla))
\end{equation}
in the derived category.

For $L = W$ (respectively ${\frak{g}}$, $CC^{per}_\bullet(W)$, $\widehat\cO$,
$\bbC$, $\widehat\Omega^\bullet$, $CC^{per}_\bullet(\widehat\cO)$),
$\ker(\nabla) = \bbA^t_M$ (respectively $Der(\bbA^t_M)$,
$CC^{per}_\bullet(\bbA^t_M)$, $\cO_M$, $\bbC_M$, $\Omega^\bullet_M$,
$CC^{per}_\bullet(\cO_M)$). We leave it to the
reader to identify $\ker(\nabla)$ for other (complexes of) ${\frak{g}}$-modules
which appear in Section \ref{section:formalRR} by analogy with the above
examples.

\subsubsection{The trace density}  \label{thetracedensity}
The image of $\tilde\mu^t$ (defined in Lemma \ref{lemma:formalTR})
under (the map on cohomology in degree zero induced by)
\begin{multline*}
GF : C^\bullet({\frak{g}},{\frak{h}};\Hom^\bullet
(CC^{per}_\bullet(W)[t^{-1}],\widehat\Omega^\bullet [2d]
[t^{-1},t]][u^{-1},u]]) @>>> \\
\R\Gamma(M;\shHom^\bullet(CC^{per}_\bullet(\bbA^t_M)[t^{-1}],
\Omega^\bullet_M [2d][t^{-1},t]][u^{-1},u]]))
\end{multline*}
is the morphism $\tilde\mu^t_{\bbA}$ defined in \cite{BNT1}. 

Similarly, the image of $\tilde\mu^t\smile\iota$ under $GF$ is the
morphism $\tilde\mu^t_{\bbA}\circ\iota$.

\subsubsection{The symbol and the Hochschild-Kostant-Rosenberg map} \label{HKR}
The image of $\sigma$ (defined in \ref{sssection:HKR})
under (the map on cohomology in degree zero induced by)
\[
GF : C^\bullet({\frak{g}},{\frak{h}};
\Hom^\bullet(CC^{per}_\bullet(W), CC^{per}_\bullet(\widehat\cO))
@>>>
\]
\[
\R\Gamma(M;\shHom^\bullet(CC^{per}_\bullet(\bbA^t_M),
CC^{per}_\bullet(\cO_M))
\]
is the morphism
\[
\sigma : CC^{per}_\bullet(\bbA^t_M) @>>> CC^{per}_\bullet(\cO_M)\ .
\]

The image of $\tilde\mu$ (defined in \ref{sssection:HKR})
under (the map on cohomology in degree zero induced by)
\begin{multline*}
GF : C^\bullet({\frak{g}},{\frak{h}};\Hom^\bullet
(CC^{per}_\bullet(\widehat\cO),\widehat\Omega^\bullet[2d][u^{-1},u]]))
@>>> \\
\R\Gamma(M;\shHom^\bullet(CC^{per}_\bullet(\cO_M),
\Omega^\bullet_M[2d][u^{-1},u]]))
\end{multline*}
is the Hochschild-Kostant-Rosenberg morphism $\tilde\mu_{\cO}$.

\subsubsection{The characteristic class of the deformation} \label{theta}
The image of the cocycle $\theta\in C^2({\frak{g}},{\frak{h}};
\frac{1}{t}\bbC[[t]])$ (defined in \ref{sssection:charcl})
under the map
\[
GF: C^\bullet({\frak{g}},{\frak{h}};\frac{1}{t}\bbC[[t]]) @>>>
\frac{1}{t}\cA^\bullet_M[[t]]
\]
is the characteristic class $\theta$ of the deformation quantization
$\bbA^t_M$ defined in \cite{F} and \cite{D} (cf. \ref{ss:theta}).

\subsubsection{The $\widehat A$-class}
The composition
\[
\widehat\Symm^\bullet({\frak{h}})^{{\frak{h}}} @>{CW}>>
C^\bullet({\frak{g}},{\frak{h}};\bbC) @>{GF}>> \cA_M^\bullet
\]
is easily seen to be the usual Chern-Weil homomorphism. In particular
we have
\[
GF(\widehat A) = \widehat A(TM)\ ,
\]
where $\widehat A$ is defined in \ref{sssection:formalA}.

Combining the above facts we obtain the following proposition.

\begin{prop}
The image under the map $\operatorname{GF}$, cf. \ref{HKR}, of the cocycle
$\tilde\mu^t\smile\iota - (\tilde\mu \circ \sigma)\smile\widehat
A\smile e^\theta$ (see Theorem
\ref{thm:formalRR}) is the morphism
$\iota\circ\tilde\mu^t_{\bbA}-(\tilde\mu_{\cO} \smile \widehat
A\smile e^\theta)\circ\sigma$.
\end{prop}

\begin{cor}\label{cor:reduction}
Theorem \ref{thm:formalRR} implies Theorem 4.6.1 of \cite{BNT1}.
\end{cor}
\section{the proof of Theorem \ref{thm:formalRR}} \label{section:formalRR-pf}
Let $W$ be the Weyl algebra of the standard symplectic space ${\bbC}^{2d}$, ${{\frak{g}}}=\text{Der}(W),$ ${{\frak{h}}}={{{\frak{sp}}}}(2d)$. Let $\tilde{\mu}^t$ be the trace density map as in \ref{tracedensity}. Note that the unit $1$ of $W$ can be viewed as a ${\frak{g}}$-invariant 
periodic cyclic zero-cycle of $W$, therefore as a zero-cocycle in $C^{\bullet}({\frak{g}}, {\frak{h}}; CC^{\operatorname{per}}_{\bullet}(W))$. The main difficulty in proving Theorem \ref{thm:formalRR} is to find $\tilde{\mu}^t(1)$
 . The idea is to replace $1$ by a homologous cycle on which $\tilde{\mu}^t$ is easy to compute. This cycle will be obtained as a result of evaluating on $1$ an operation acting on $C^{\bullet}({{\frak{g}}}, {{\frak{h}}}; CC^{\operatorname{per}}_{\bullet}(W[t^{-1}]))$. Complexes of such operations are studied in 
\ref{ss:opera}. In particular, it is shown in Theorem \ref{thm:calc3} that for any algebra $A$ the reduced cyclic complex $(\overline{C}^{\lambda}_{\bullet}(A))[[u,u^{-1}], ub) $ acts on $CC^{\operatorname{per}}_{\bullet}(A)$. Motivated by this, we construct the fundamental class $U$ in the hypercohomology $H^{\bullet}({{\frak{g}}}, {{\frak{h}}}; \overline{C}^{\lambda}_{\bullet}(W))[t ^{-1}][[u,u^{-1}] $ of total degree $1$ (subsection \ref{ss:fundclass}).  We will show that $\tilde{\mu} ^{t} \cdot (U \bullet ?) = \text{Id}$ where $\bullet$ is the pairing from Theorem \ref{thm:calc3} (Lemma \ref{lemma:muhumu1}). Therefore, to prove Theorem \ref{thm:formalRR} we have to show that the operation $U\bullet ?$ is cohomologous to $(\widehat{A}\cdot e^{\theta})^{-1} \smile ?$. To do that, in subsection \ref{ss:computu} we compute the image of the class $U$ in $H^{\bullet}({{\frak{g}}}, {{\frak{h}}}; \overline{C}^{\lambda}_{\bullet}(W[\eta]))[t ^{-1}][u,u^{-1}] $ in terms of classes $\eta ^{(m)}$ (cf. subsections \ref{ss:ccdga}, \ref{ss:boma}). 

All the facts about Hochschild, cyclic, and Lie algebra (co)homologies that we use is contained in sections \ref{section:CC} and \ref{section:lie}.

\subsection{The fundamental class} \label{ss:fundclass}
\begin{lemma} \label{lemma:ccredw}
Let $W$ denote the Weyl algebra of the symplectic vector space ${\Bbb{C}} ^{2d}$. Then
$$\overline{HC}_i(W[t^{-1}]) = {\Bbb{C}}[[t,  t^{-1}], \;\; i = 1, 3, \cdots, 2d-1;
$$
$$\overline{HC}_i(W[t^{-1}]) =  0
$$
otherwise. (\cite{Bry}, \cite{FT1}).
\end{lemma}
\begin{pf} From proposition \ref{prop:HH(W)} and from the spectral sequence (\ref{eq:ssccc}) it follows immediately that 
\begin{equation} \label{eq:HC(W)}
{HC}_i(W[t^{-1}]) = {\Bbb{C}}[[t,  t^{-1}],\; i = 2d+2i, \;i\geq 0
\end{equation}
\begin{equation} \label{eq:HC(W)1}
{HC}_i(W[t^{-1}]) = 0
\end{equation}
otherwise. The explicit formula for a generator of $HH_{2d}(W[t^{-1}])$ (after (\ref{map:K-to-DR})) implies that the image of this generator in $HC_{2d}(W[t^{-1}])$ goes to zero under the connecting morphism $\partial$ to the reduced cyclic homology. Therefore the map
\begin{equation} \label{HC(k)-to-HC(W)}
{HC}_i( {\Bbb{C}}[[t,  t^{-1}]) \to {HC}_i(W[t^{-1}])
\end{equation}
is an isomorphism for $i=2d$. Since this map commutes with the Bott periodicity map $S$, one sees that it is an isomorphism for all $i \geq 2d$. The lemma follows now from the long exact sequence associated to the exact triangle (\ref{ses:CC-CCred}).
\end{pf}
Let
\begin{equation} \label{ses:u0}
U_0 = \frac{u^{-d}}{2d\cdot t^d}\text{Alt}(v_1 \otimes \cdots \otimes v_{2d})
\end{equation}
where $\{v_1, \cdots, v_{2d}\}$ is a Darboux basis of $V^*.$ A simple computation shows that $\text{Br}(U_0) = 1$ where  $\text{Br}$ is defined in (\ref{ses:jacek}). 
\begin{lemma} \label{lemma:funclass}
The cycle $U_0$ extends to a cocycle $U$ of the complex $C^{\bullet}({{\frak{g}}}, {{\frak{h}}}; \overline{C}^{\lambda}_{\bullet}(W))[t ^{-1}][u,u^{-1}] $ of total degree $1-2d.$ Any two such extensions are cohomologous.
\end{lemma}
\begin{pf}
Follows immediately from the fact that $\overline{HC}_i(W[t^{-1}]) =  0$ for $i>2d-1.$
\end{pf}
It follows that the cohomology class of $U$ is canonically defined. We will call this cohomology class {\it {the fundamental class}}.

\subsection{Computation of the fundamental class} \label{ss:computu}
For any algebra $A$, consider the DGA
$A[\eta]\overset{def}{=}
(A,\delta)\otimes_k 
(k[\eta], \displaystyle\frac{\partial}{\partial \eta})$, where
$\deg\eta = 1$ (in particular $\eta^2=0$) (cf. \ref{ss:ccdga}). Let
$$
\eta ^{(n+1)} \overset{def}{=} n!\eta ^{\otimes n+1}
$$
be the cycle in the reduced cyclic complex $\overline{C}^{\lambda}_{2n+1}(A[\eta])$, as in \ref{ss:boma}.

\begin{thm} \label{thm:prou1}
One has 
$$
U = \sum_{m\geq 0} (\widehat{A} \cdot e^{\theta})^{-1}_{2m} \cdot \eta ^{(m+d)} \cdot u^{-m}
$$
in  $H^{\bullet}({{\frak{g}}}, {{\frak{h}}}; \overline{C}^{\lambda}_{\bullet}(W[\eta]))[t ^{-1}][u,u^{-1}]. $
\end{thm}
\begin{pf}
Recall that $\tilde{{\frak{g}}}=\frac{1}{t}W$ is a central extension of ${{\frak{g}}} = \text{Der}(W).$ The preimage of ${{\frak{h}}} = {{\frak{sp}}}({2d})$ in $\tilde{{\frak{h}}}$ splits canonically:  $\tilde{{\frak{h}}}={{\frak{h}}}\oplus  \frac{1}{t}{\Bbb{C}}[[t]].$ 

(The subalgebra ${{\frak{h}}}$ inside $\tilde{{\frak{h}}}$ is identified with the linear span of the elements $\frac{1}{t}v \ast w + \frac{1}{2}\omega(v,w)$ where $v$ and $w$ are in $V^{\ast}$). Note that for any $({{\frak{g}}},{{\frak{h}}} )$-module $L$
$$ C^{\bullet}( {{\frak{g}}},{{\frak{h}}}; L) = C^{\bullet}( \tilde{{\frak{g}}}, \tilde{{\frak{h}}}; L).
$$
Note that the complexes of ${{\frak{g}}}$-modules
$ \overline{C}^{\lambda}_{\bullet}(W))[t ^{-1}][u,u^{-1}]$ and
$\overline{C}^{\lambda}_{\bullet}(W[\eta]))[t ^{-1}][u,u^{-1}] $ are homotopically constant in sense of Definition \ref{hotri}. Indeed, for any $\Phi \in W[\eta]$ and for $a_i \in A[\eta]$ let us define$$\iota_{\Phi}(a_0 \otimes \cdots \otimes a_p)    =  $$
\begin{equation}  \label{ses:ifi}
=\frac{u^{-1}}{t}\sum_{i=0}^{p}(-1)^{\sum_{p \leq i}(\text{deg}a_p+ 1)( \text{deg}\Phi + 1)} 
 a_0 \otimes \cdots a_i \otimes \Phi \otimes \ldots \otimes a_p
\end{equation}
 if $X=\frac{1}{t}\text{ad}(\Phi)$ then one can put $\iota_X = \iota_{\Phi}.$

We would like to apply the generalized Chern-Weil map from (\ref{ss:charactlie}) and reduce the theorem to a computation in $C^{\bullet}( \tilde{{\frak{h}}}[\epsilon], \tilde{{\frak{h}}}; \overline{C}^{\lambda}_{\bullet}(W))[t ^{-1}][u,u^{-1}])$.
The problem is, the fundamental class $U$ can be easily defined in this context but the cochains $\eta ^{(m)}$ are not ${{\frak{h}}}[\epsilon]$-invariant.To correct this, let us introduce cocycles in $C^{\bullet}(\tilde{\frak{g}}[\epsilon], \tilde{\frak{h}}; \overline{C}_{\bullet}^{\lambda}(W[\eta])[[u, u^{-1}])$:
\begin{equation} \label{ses:neweta}
(X_1\epsilon, \ldots, X_p\epsilon) \mapsto \iota_{X_1 \eta} \ldots \iota_{ X_p \eta} \eta ^{(m)}
\end{equation}
for all $p \geq 0$;
\begin{equation} \label{ses:neweta1}
(X_1, \ldots ) \mapsto 0, \;X_i \in {{\frak{g}}}.
\end{equation}
Note that the restriction from $C^{\bullet}(\tilde{{\frak{g}}}[\epsilon], \tilde{{\frak{h}}}; ...)$ to $C^{\bullet}(\tilde{{\frak{h}}}[\epsilon], \tilde{{\frak{h}}}; ...)$ is a quasi-isomorphism. We will denote by $\eta ^{[m]}$ the image of the cocycle (\ref{ses:neweta}, \ref{ses:neweta1}) under this restriction. Let us observe that there are two cocycles in $C^{\bullet}(\tilde{{\frak{g}}}[\epsilon], \tilde{{\frak{h}}};\ldots)$ whose restriction to $C^{\bullet}(\tilde{{\frak{h}}}[\epsilon], \tilde{{\frak{h}}};\ldots)$ is $\eta ^{[m]}:$ one is ${\bf {CW}}(\eta ^{[m]})$, the other is the cocycle (\ref{ses:neweta}, \ref{ses:neweta1}). But the restrictions of those cocycles to $C^{\bullet}(\tilde{{\frak{g}}}, \tilde{{\frak{h}}};\ldots)$ are respectively $CW (\eta ^{[m]})$ and $\eta ^{(m)}.$ We therefore see that 
\begin{equation} \label{ses:etas}
CW (\eta ^{[m]}) = \eta ^{(m)}
\end{equation}
in $H^{\bullet}( \tilde{{\frak{g}}}, \tilde{{\frak{h}}}; \overline{C}^{\lambda}_{\bullet}(W))[t ^{-1}][u,u^{-1}]$. Thus, it suffices to prove the following
\begin{lemma} \label{lemma:prouvh}
One has 
$$
U = \sum_{m\geq 0} (\widehat{A} \cdot e^{\theta})^{-1}_{2m} \cdot \eta ^{[m+d]} \cdot u^{-m}
$$
in  $H^{\bullet}(\tilde{{\frak{h}}}[\epsilon], \tilde{{\frak{h}}}; \overline{C}^{\lambda}_{\bullet}(W[\eta]))[t ^{-1}][[u,u^{-1}]. $
\end{lemma}
To prove this, let us first reduce the statement to the case $d=1.$ Let ${{\frak{d}}}_d$ be the subalgebra of diagonal matrices in ${{\frak{sp}}}({2d})$; 
$\tilde{{\frak{d}}}_d = {{\frak{d}}}_d \oplus \frac{1}{t}{\Bbb{C}}[[t]].$
\begin{lemma} \label{lemma:splitting} (the splitting principle).
The restriction $H^{\bullet}(\tilde{{\frak{h}}}[\epsilon], \tilde{{\frak{h}}}; k) @>>>H^{\bullet}(\tilde{{\frak{d}}}_d[\epsilon], \tilde{{\frak{d}}}_d; k)$ is a monomorphism.
\end{lemma}
\begin{pf}
Indeed, the left hand side is the space of ${\frak{h}}$-invariant polynomials on ${\frak{h}}^{\ast}$ tensored by the symmetric algebra of the space $\frac{1}{t}\C [[t]]$ whereas the right hand side is the space of polynomials on ${\frak{d}}^{\ast}$ tensored by the symmetric algebra of the space $\frac{1}{t}\C [[t]]$. But the space of ${\frak{h}}$-invariant polynomials on ${\frak{h}}^{\ast}$ maps to the space of polynomials on ${\frak{d}}^{\ast}$ as the subspace of polynomials invariant under the Weyl group.
\end{pf}

Therefore it is enough to prove Lemma \ref{lemma:prouvh} if one replaces $\tilde{{\frak{h}}}$ by $\tilde{{\frak{d}}}_d.$ Indeed, we want to compute $U$ in the cohomology which is isomophic to $H^{\bullet}(\tilde{{\frak{h}}}[\epsilon], \tilde{{\frak{h}}}; CC_{\bullet}(k));$ but $CC_{\bullet}(k) \simeq k[u^{-1},u]]/k[[u]]$, thus the map on the cohomology which is induced by $\tilde{\frak{d}}_d \rightarrow \tilde{\frak{h}}$ is a monomorphism.

Next, note that ${{\frak{d}}}_d = \oplus ^d {{\frak{d}}}_1$ and $W=\otimes ^d W_1$ where $W_1$ is the Weyl algebra of the standard two-dimensional symplectic space.
Now, note that $U=\times {^d} U_1$ where $U_1$ is the fundamental class for $d=1$ and $\times$ is the external product defined in \ref{sss:cross}. Note also that $\eta^{[m]} = \eta ^{\times m}$ for any $m.$ So it is enough to prove lemma \ref{lemma:prouvh} if one replaces $d$ by $1$ and ${{\frak{h}}}$ by ${{\frak{d}}}_1.$

But in this case the fundamental class $U$ can be constructed explicitly. Put $\partial = \frac{1}{t} \xi;$ then $U$ is represented by the cocycle 
\begin{equation} \label{ses:uin1}
U = \sum_{m=1}^{\infty}\frac{u^{1-2m}}{m}(\partial \otimes x)^{\otimes m}c_1^{m-1}
\end{equation}
Recall that $\tilde{{\frak{d}}}_1 = {\Bbb{C}} \oplus \frac{1}{t} {\Bbb{C}}[[t]].$ For $(a,b) \in {\Bbb{C}} \oplus \frac{1}{t} {\Bbb{C}}[[t]]$ let 
$$ c_1(a,b) = a,/;; \theta (a,b) = b$$
 One has 
$$ C^{\bullet}( \tilde{{\frak{d}}}_1 [\epsilon], \tilde{{\frak{d}}}_1; \overline{C}^{\lambda}_{\bullet}(W[\eta]))[t ^{-1}][[u,u^{-1}] \simeq \overline{C}^{\lambda}_{\bullet}(W[\eta]))[t ^{-1}][[u,u^{-1}][[c_1, \theta]]
$$
with the differential $\frac{\partial}{\partial \eta} + ub + c_1 \cdot \iota_{x * \partial}.$ The operator $\iota_{x * \partial}$ is given by formula (\ref{ses:ifi}).
We would like to compute $U$ as in subsection \ref{ss:boma}. For this we can use
Lemma \ref{lemma:ch-jacek}. Take $j: W@>>> {\Bbb{C}}[[t,t^{-1}],$ $j(f(x,\xi,t))=f(0,0,0).$
One has $\rho(\partial, x)=\frac{1}{2};\;\;\rho(x, \partial)=-\frac{1}{2}.$ Note that $\text{Br} = 0$ on the image of $\iota_{x * \partial}.$ Indeed, $\rho({x * \partial},a)-\rho(a,{x * \partial})=0$ for all $a.$ Therefore the map $\text{Br}$ is a morphism of complexes.
One has
\begin{eqnarray*}
\text{Br}(U)=\sum_{m=1}^{\infty}\frac{u^{1-m}}{m}(( \frac{1}{2})^m - (-\frac{1}{2})^m)\frac{m}{m!}1^{(m)}c_1^{m-1} = \\
\sum_{m=1}^{\infty}u^{1-m}1^{(m)}[\frac{e^{\frac{c_1}{2}}-e^{-\frac{c_1}{2}}}{c_1}]_{2m-2}
\end{eqnarray*}
On the other hand, 
$$\text{Br}(\eta^{[m+1]}) = \sum_{l=0}^{\infty}u^{-l}1^{(m+l)}(e^{\theta})_{2l}$$
for every $m.$ Therefore in the cohomology
\begin{equation}
U=\sum_{m=1}^{\infty}(\frac{e^{\frac{c_1}{2}}-e^{-\frac{c_1}{2}}}{c_1}e^{-\theta})_{2m}u^{-m}\eta^{[m+d]}
\end{equation}
because their images under the map $\text{Br}$ are the same and because this map is a quasi-isomorphism by lemma \ref{lemma:ccredeta}. This proves Theorem \ref{thm:prou1}.
\subsection{End of the proof of Theorem \ref{thm:formalRR}}.
\begin{lemma} \label{lemma:muhumu}
If one regards $1$ in the left hand side as an element of $C^0({{\frak{g}}},{{\frak{h}}}; CC^{\text{per}}_{0}(W))$, then 
$$ \tilde{\mu}^{t} (U\bullet 1) = 1$$
in $C^{\bullet}({{\frak{g}}},{{\frak{h}}}; \Omega^{\bullet}[2d][t^{-1},t]][u^{-1},u]]))$
\end{lemma}
\begin{pf} Because of Theorem \ref{thm:calc3}, 4, of formula (\ref{ses:u0}), and of the definition of $U$, one can conclude the following. The cochain $U\bullet 1$ has the component in $C^0 ( {{\frak{g}}},{{\frak{h}}}; u^{-d}C_{2d})$ which is equal to 
$$(U\bullet 1)_0 = \frac{u^{-d}}{t^d}1 \otimes \operatorname{Alt}(v_1 \otimes \ldots \otimes v_{2d})$$
where $v_1, \ldots, v_{2d}$ is a Darboux basis; all other components of $U \bullet 1$ are in $C^p ({{\frak{g}}},{{\frak{h}}}; u^{i}C_{2d+j})$ where $p>0$ and $j>0$. Therefore these other components are sent to zero by $\tilde{\mu}^t$, whereas 
$$\tilde{\mu}^t((U\bullet 1)_0)=1$$
\end{pf}
\begin{lemma} \label{lemma:muhumu1}
If one regards $1$ in the left hand side as an element of $H^0({{\frak{g}}},{{\frak{h}}}; CC^{\text{per}}_{0}(W))$, then 
$$ \tilde{\mu}^{t} (1) = \sum (\widehat{A} \cdot {e^{\theta}})_{2m}u^{d-m}$$
in $C^{\bullet}({{\frak{g}}},{{\frak{h}}}; \Omega^{\bullet}[2d][t^{-1},t]][u^{-1},u]]))$
\end{lemma}
\begin{pf}
Follows immediately from Theorem \ref{thm:prou1} and Theorem \ref{thm:calc3}, 1.
\end{pf}
Now theorem \ref{thm:formalRR} follows from the above lemma. Indeed, the inclusion of $\C[u^{-1},u]]$ into $CC^{\text{per}}_{\bullet}(W))$ is a $\frak{g}$-equivariant quasi-isomorphism because of Theorem \ref{thm:rigid}. This inclusion induces a quasi-isomorphism 
$$C^{\bullet}({{\frak{g}}},{{\frak{h}}}; \operatorname{Hom}^{\bullet}(CC^{\text{per}}_{\bullet}(W), \Omega^{\bullet}[2d][t^{-1},t]][u^{-1},u]])) \longrightarrow$$
$$C^{\bullet}({{\frak{g}}},{{\frak{h}}}; \Omega^{\bullet}[2d][t^{-1},t]][u^{-1},u]]))$$
(recall that all our homomorphisms are by definition $\C[u^{-1},u]]$-linear and $(u)$-adically continuous).
But under this quasi-isomorphism the images of the two sides of Theorem \ref{thm:formalRR} are the same, namely the left hand side of the formula in Lemma \ref{lemma:muhumu1}.
\end{pf}

\section{Appendix: Review of Hochschild and cyclic homology}
\label{section:CC}
In this section we review the basic definitions of the Hochschild and the
(negative, periodic, reduced) cyclic complex of an algebra, of a sheaf of algebras
on a space as well as of the ``topological'' versions of the above. In
particular we establish notational conventions with regard to Hochschild and
cyclic complexes which are used in the body of the paper. After that, we will review additional results and constructions that we used in section \ref{section:formalRR-pf}, such as: an explicit formula for the connecting morphism from the reduced cyclic complex, an external multiplication on the reduced cyclic complex, and operations on Hochschild and cyclic complexes. 
\subsection{Hochschild and cyclic complexes of algebras}
\label{subsection:C-CC}
Let $k$ denote a commutative algebra over a field of characteristic zero
and let $A$ be a flat $k$-algebra with $1_A\cdot k$ contained in the center,
not necessarily commutative. Let $\overline A = A/k$, and let
$C_p(A)\overset{def}{=}A\otimes_k\overline A^{\otimes_k p}$. Define
\begin{eqnarray}\label{diffl:b}
b : C_p(A) & @>>> & C_{p-1}(A) \\ \nonumber
a_0\otimes\cdots\otimes a_p & \mapsto & 
(-1)^p a_pa_0\otimes\cdots\otimes a_{p-1} + \\ & &
\sum_{i=0}^{p-1}(-1)^ia_0\otimes\cdots\otimes a_ia_{i+1}\otimes\cdots\otimes
a_p\ .\nonumber
\end{eqnarray}
Then $b^2 = 0$ and the complex $(C_\bullet, b)$, called {\em the standard
Hochschild complex of $A$} represents $A\otimes^{\bbL}_{A\otimes_kA^{op}}A$
in the derived category of $k$-modules. The homology of this complex is denoted by $H_{\bullet}(A,A)$, or by $HH_{\bullet}(A)$.

The map
\begin{eqnarray}\label{diffl:B}
B : C_p(A) & @>>> & C_{p+1}(A) \\ \nonumber
a_0\otimes\cdots\otimes a_p & \mapsto & \sum_{i=0}^p (-1)^{pi}
1\otimes a_i\otimes\cdots\otimes a_p\otimes a_0\otimes\cdots\otimes a_{i-1}
\end{eqnarray}
satisfies $B^2 = 0$ and $[B,b] = 0$ and therefore defines a map of complexes
\[
B: C_\bullet(A) @>>> C_\bullet(A)[-1]\ .
\]
For $i,j,p\in\bbZ$ let
\begin{eqnarray*}
CC^-_p(A) & = & \prod_{\overset{i\geq p}{i = p \mod 2}} C_{i}(A) \\
CC^{per}_p(A) & = & \prod_{i = p \mod 2} C_{i}(A)\ .
\end{eqnarray*}
$$CC_p(A)  =  \bigoplus_{\overset{i\leq p}{i = p \mod 2}} C_{i}(A)$$
The complex $(CC^-_\bullet(A),B+b)$ (respectively
$(CC^{per}_\bullet(A),B+b)$, respectively $(CC_\bullet(A),B+b)$) is called the {\em negative cyclic} (respectively
{\em periodic cyclic}, respectively {\em cyclic}) complex of $A$. The homology of these complexes is denoted by $HC^-_\bullet(A)$, respectively by $HC^{per}_\bullet(A)$, respectively by $HC_\bullet(A)$.

There are inclusions of complexes
\begin{equation}\label{inclusions}
CC^-_\bullet(A)[-2]\hookrightarrow CC^-_\bullet(A)\hookrightarrow
CC^{per}_\bullet(A)
\end{equation}
and the short exact sequences
\begin{equation}\label{ses:CC-C}
0 @>>> CC^-_\bullet(A)[-2] @>>> CC^-_\bullet(A) @>>> C_\bullet(A) @>>> 0\ .
\end{equation}
\begin{equation}\label{ses:CC-C1}
0 @>>>  C_\bullet(A) @>>> CC_\bullet(A) @>S>> CC_\bullet(A) [2] @>>> 0\ .
\end{equation}
The operator $S$ is called the Bott periodicity operator. To the double complex $CC_{\bullet}(A)$ one associates the spectral sequence
\begin{equation} \label{eq:ssccc}
E^2_{pq} = H_{p-q}(A, A)
\end{equation}
converging to $HC_{p+q}(A)$.

In what follows we will use the notation of Getzler and Jones (\cite{GJ}).
Let $u$ denote a variable of degree $-2$ (with respect to the homological
grading). Then the negative and periodic
cyclic complexes are described by the following formulas:
\begin{eqnarray} \label{ses:gejo}
CC^-_{\bullet}(A) & = & (C_{\bullet}(A)[[u]], b + uB) \\
\label{ses:gejoper}
CC^{\text{per}}_{\bullet}(A) & = & (C_{\bullet}(A)[[u, u^{-1}], b + uB)\ .
\end{eqnarray}
\begin{equation} \label{eq:ccu}
CC_{\bullet}(A)  =  (C_{\bullet}(A)[[u, u^{-1}] / C_{\bullet}(A)[[u]], b + uB)
\end{equation}
In this language, the Bott periodicity map $S$ is just multiplication by $u$.
\subsection{The reduced cyclic complex} \label{ss:recyco}
Recall the original definition of the cyclic complex from \cite{C1}, \cite{Ts1}. Let
$\tau = \tau_p$ denote the endomorphism of $A^{\otimes_kp+1}$
given by the formula 
\[
\tau (a_0\otimes\cdots\otimes a_p) =
(-1)^p a_p\otimes a_0 \cdots \otimes a_{p-1}
\]
Let
\begin{equation}
C^{\lambda}_p(A) = A^{\otimes_kp+1}/{\operatorname{Im}}(\id - \tau)\ .
\end{equation}
 Recall that, if one defines 
\begin{eqnarray*}
b' : C_p(A) & @>>> & C_{p-1}(A) \\
a_0\otimes\cdots\otimes a_p & \mapsto & 
\sum_{i=0}^{p-1}(-1)^ia_0\otimes\cdots\otimes a_ia_{i+1}\otimes\cdots\otimes
a_p\ 
\end{eqnarray*}
and 
\begin{eqnarray*}
N : C_p(A) & @>>> & C_{p-1}(A) \\
N = \operatorname{id} + \tau + \ldots \tau ^{p-1},
\end{eqnarray*}
then
\begin{equation} \label{ses:ntb}
b (\operatorname{id} - \tau) = (\operatorname{id} - \tau)b'; \;\; b'N=Nb.
\end{equation}
Therefore 
the differential $b$ descends to a map
\[
b: C^{\lambda}_p(A)  @>>> C^{\lambda}_{p-1}(A)\ .
\]

The complex $C^{\lambda}_\bullet (A)$ is isomorphic to $CC_{\bullet}(A)$ in the derived category. (Sketch of the proof: in the cyclic double complex $CC_{\bullet}(A)$ one can replace $C_{\bullet}(A)$ by the cone of $\operatorname{id}-\tau$ and $B$ by $N$. One gets a double complex which is quasi-isomorphic to $CC_{\bullet}$ and which is in turn quasi-isomorphic to $C^{\lambda}_{\bullet}$. Cf. \cite{L} for more detail). 
Let $\overline{A} = A{/}k$ and let
\begin{equation}\label{ses:cyclicred}
\overline{C}^{\lambda}_p (A) =
\overline{A}^{\otimes p+1}/{\operatorname{Im}}(\id - \tau)\ .
\end{equation}
It is easy to see that the diferential $b$ descends to
$\overline{C}^{\lambda}_\bullet (A)$. We denote the homology of the complex
$\overline{C}^{\lambda}_\bullet (A)$ by $\overline{HC}_{\bullet}(A)$.
There is an exact triangle
\begin{equation}\label{ses:CC-CCred}
C^{\lambda}_{\bullet} (k) @>>> C^{\lambda}_{\bullet} (A) @>>>
\overline{C}^{\lambda}_\bullet (A) @>>> C^{\lambda}_{\bullet} (k)[1]\ .
\end{equation}
This follows from the Hochschild-Serre spectral sequence 
$$E^2_{pq} = H_p ({\frak{gl}}(A), {\frak{gl}}(k)) \otimes H_q ({\frak{gl}}(k))$$
converging to $H_{p+q}({\frak{gl}}(A))$ and from theorem \ref{thm:liecyclic}.

\subsection{Homology of differential graded algebras}  \label{ss:ccdga}
One can easily generalize all the above constructions to the case when $A$ is
a differential graded algebra (DGA) with the differential $\delta$ (i.e
$A$ is a graded algebra and $\delta$ is a derivation of degree $1$ such that
$\delta^2=0$). 

The action of $\delta$ extends to an action on Hochschild chains by the
Leibnitz rule:
\[
\delta (a_0\otimes\cdots\otimes a_p ) = 
\sum_{i=1}^p {(-1)^{\sum_{k<i}{(\deg a_k + 1)+1}}
(a_0\otimes\cdots\otimes\delta a_i \otimes \cdots \otimes a_p )}
\]
The maps $b$, $B$, $\tau$ are modified to include signs (cf. \cite{NT2}).
The complex $C_{\bullet}(A)$ (respectively $C^{\lambda}_{\bullet}(A)$,
$\overline{C}^{\lambda}_{\bullet}(A)$) now becomes the total complex of
the double complex with the differential $b + \delta$. The negative and
the periodic cyclic complexes are defined as before in terms of the new
definition of $C_{\bullet}(A)$.

Given a DGA $(A,\delta)$, consider the DGA
$A[\eta]\overset{def}{=}
(A,\delta)\otimes_k 
(k[\eta], \displaystyle\frac{\partial}{\partial \eta})$, where
$\deg\eta = 1$ (in particular $\eta^2=0$).
 
\begin{lemma} \label{lemma:ccredeta}
The complex  $C^{\lambda}_{\bullet}(A[\eta])$ is acyclic.
\end{lemma}
\begin{cor}\label{cor:bdry-map}
The morphism
\[
\overline{C}^{\lambda}_\bullet (A[\eta]) @>>> C^{\lambda}_{\bullet} (k)[1]
\]
is an isomorphism.
\end{cor}

\subsection{Relation to Lie algebra homology} \label{ss:liecyclic}

For a DGA $(A,\delta)$ over $k$ let 
${\frak{gl}}(A) = \underset{n}{\varinjlim}{\frak{gl}}_n (A)$ denote the
(DGLA) of finite matrices. Note that ${\frak{gl}}(k)$ is a DGL-subalgebra
of ${\frak{gl}}(A)$. Theorem \ref{thm:liecyclic} below
identifies their respective subcomplexes of primitive elements of the DG
coalgebras $C_{\bullet}({\frak{gl}}(A))_{{\frak{gl}}(k)}$
and $C_{\bullet}({\frak{gl}}(A),{\frak{gl}}(k))$.

Let $E^a_{pq}$ denote the elementary matrix with
$\left(E^a_{pq}\right)_{pq}=a$ and other entries equal to zero.

\begin{thm} \label{thm:liecyclic}
The map
\begin{eqnarray*}
A^{\otimes p+1} & @>>> & \bigwedge^{p+1}{\frak{gl}}(A) \\
a_0\otimes\cdots\otimes a_p & \mapsto &
E_{01}^{a_0} \wedge  E_{12}^{a_1} \wedge\dots\wedge E_{p-1, 0}^{a_p}
\end{eqnarray*}
induces isomorphisms of complexes
\begin{eqnarray*}
C^{\lambda}_\bullet(A) & @>>> & 
\operatorname{Prim} C_{\bullet}({\frak {gl}}(A))_{{\frak {gl}}(k)}[1] \\
\overline{C}^{\lambda}_\bullet(A) & @>>> &
\operatorname{Prim} C_{\bullet}({\frak {gl}}(A), {\frak {gl}}(k))[1]
\end{eqnarray*}
\end{thm}
The proof is contained in \cite{Ts}, \cite{LQ}, \cite{L}

\subsection{The connecting morphism} \label{ss:boma}
We will now give an explicit construction of the connecting morphism
in the exact triangle \eqref{ses:CC-CCred} and the isomorphism (in the
derived category) of Corollary \ref{cor:bdry-map}. Because of the application that we need (via Theorem \ref{thm:calc3}), we will work with a slightly generalized situation. Namely, let $u$ be a formal parameter of degree $2$. For any $DGA (A, \; \delta)$, we consider the complex $(\overline{C}^{\lambda}_\bullet(A)[u], \; ub+\delta)$. Define the boundary map $\partial$ as the composition
\begin{equation} \label{eq:partialu}
(\overline{C}^{\lambda}_\bullet(A)[u], \; ub+\delta) \longrightarrow (\overline{C}^{\lambda}_\bullet(A[\eta])[u], \; ub+\delta + \frac{\partial}{\partial \eta}) \longrightarrow \end{equation}
$$
\longrightarrow (\overline{C}^{\lambda}_\bullet(k[\eta])[u], \; \frac{\partial}{\partial \eta}) \stackrel{\backsim}{\longrightarrow} C^{\lambda}_{\bullet}(k)[1][u]
$$
where the second morphism is an inverse to the quasi-isomorphism induced by the embedding $k\rightarrow A$ and the last one is the boundary map in the triangle (\ref{ses:CC-CCred}) (the composition of the last two is as in Corollary \ref{cor:bdry-map}).

It is easy to see that the map from Theorem \ref{thm:liecyclic} induces an isomorphism
\begin{equation} \label{eq:LQTu}
({\overline{C}}^{\lambda}_\bullet(A),ub +\delta)  @>>> 
\operatorname{Prim} C_{\bullet}({\frak {gl}}(A), {\frak {gl}}(k))[1], u\partial _{\operatorname{Lie}} + \delta)
\end{equation}

For a DGA $(A,\delta)$ we denote
by $T(A)$ the tensor algebra of $A$ and by $j$ a $k$-linear map $A @>>> k$
such that $j(1)=1$.

Let $\rho : T(A) @>>> k[u]$ denote the map defined by
\begin{eqnarray*}
\rho (a) & = & \text{$j(\delta a)$ for $a \in A$,} \\
\rho (a_1 \otimes a_2) & = & \text{$u j( a_1)j(a_2) - j(a_1 a_2)$ for
$a_1,a_2\in A$, and} \\ 
\rho & = & \text{$0$ on $A^{\otimes p}$ for $p \neq 1,2$.} 
\end{eqnarray*}
For $a_i \in A$,
$a_0\otimes\cdots\otimes a_p \in A^{\otimes p+1}$ with
$\left(\sum_i(\deg a_i+1)\right)-1 = 2n+1$
let
\begin{multline} \label{ses:jacek}
\operatorname{br}_{2n+1}(a_0 \otimes \cdots \otimes a_p) = \\
\frac{1}{n!}\sum_{i=0}^{p}
(-1)^{{\sum_{k<i}(\deg a_k +1)\sum_{k\geq i}(\deg a_k +1)}} 
(\rho\otimes\dots\otimes\rho)
(a_i \otimes \cdots a_0\otimes\dots\otimes a_{i-1})\ .
\end{multline}
This defines a map
$\operatorname{br}_{2n+1}: \overline{C}^{\lambda}_{2n+1}(A)[u] @>>> k[u]$ (cf. \cite{Br}).
Let $\operatorname{br}_{2n} = 0$ on $\overline{C}^{\lambda}_{2n}(A)$.

Let $\One^{(n+1)} = n!(n+1)!\cdot 1^{\otimes 2n+1}\in C^{\lambda}_{2n}(k)$.
Let $\operatorname{Br}_{2n+1}^A=\operatorname{br}_{2n+1}\cdot\One^{(n)}$,
$\operatorname{Br}_{2n}^A=0$,
$\operatorname{Br}^A=\sum_p\operatorname{Br}_{p}$. Occasionally, when this
causes no confusion, we will omit the reference to the algebra $A$ from
the notation and simply write $\operatorname{Br}$.

This defines the map of graded $k[u]$-modules
\[
\operatorname{Br}^A : \overline{C}^{\lambda}_\bullet(A)[u] @>>>
C^{\lambda}_\bullet(k)[1][u]
\]
which is, in fact, a map of complexes as follows from the alternative
description given below.

The splitting $j$ induces the splitting
${\pmb{\rho}} : {\frak{gl}}(A) @>>> {\frak{gl}}(k)$ of the inclusion
${\frak{gl}}(k)\hookrightarrow{\frak{gl}}(A)$. Let $P_n$ denote the invariant
polynomial $X\mapsto\displaystyle\frac1{n!}tr(X^n)$ on ${\frak{gl}}(k)$.
Put 
\begin{equation}
\label{eq:Rwithu}
R({\pmb{\rho}})=(u\partial ^{\operatorname {Lie}}+\delta)R({\pmb{\rho}}) + \frac{1}{2}[R({\pmb{\rho}}),R({\pmb{\rho}})]
\end{equation}
Then $\operatorname{br}_{2n+1}$ is the Chern character cochain $P_n (R({\pmb{\rho}}))$ composed with the morphism (\ref{eq:LQTu}).

Cf. \ref{chclasstrivco} for more detail on the generalized Chern-Weil map for relative Lie algebra cohomology.

\begin{lemma} \label{lemma:ch-jacek} The restriction of the pairing with $P_n (R({\pmb{\rho}}))$ to the subspace
of all primitve elements (see Theorem \ref{thm:liecyclic}) coincides with the
map $\operatorname{br}_{2n+1}$.
\end{lemma}
The proof is straightforward. 
\begin{cor}
The map $\operatorname{Br}^A$ is a morphism of complexes.
\end{cor}

\begin{prop}
The morphism 
$\operatorname{Br}^A : \overline{C}^{\lambda}_\bullet(A) @>>>
C^{\lambda}_\bullet(k)[1]$ represents the morphism $\partial$ (\ref{eq:partialu}).
\end{prop}
\begin{pf}

\begin{definition} \label{def:etan}
Put
$$\eta^{(n+1)}\overset{def}{=}n!\cdot\eta^{\otimes n+1}
\in \overline{C}^{\lambda}_{2n+1} (k[\eta])$$
\end{definition}

\begin{lemma} \label{lemma:ooo}

\begin{itemize}
\item
The image of  $\eta^{(n+1)}$ under the connecting morphism
is equal to the class of of the cycle $\One^{(n+1)}$
\item
$\operatorname{br}_{2n+1}(\eta^{(n+1)}) = 1$
\end{itemize}
\end{lemma}
This is easily checked by a direct computation. 

The map $\partial$ is functorial, and $\operatorname{Br}^{A[\eta]}$ restricts to
$\operatorname{Br}^A$  (respectively $\operatorname{Br}^{k[\eta]}$).
Therefore, it is sufficient to verify that $\operatorname{Br}^{k[\eta]}$
represents the connecting morphism from (\ref{ses:CC-CCred}). But this follows from Lemma \ref{lemma:ooo}.

\end{pf}
It is easy to see that the specialization of the morphism (\ref{eq:partialu}) at $u=1$ is equal to the connecting morphism in the triangle \eqref{ses:CC-CCred}.
\subsection{External product on the reduced cyclic complex}
\label{sss:cross}

Here we recall the product
\begin{equation} \label{ses:prodccred}
\times: \overline{C}^{\lambda}_{\bullet}(A) \otimes \overline{C}^{\lambda}_{\bullet}(B) @>>> \overline{C}^{\lambda}_{\bullet + 1}(A \otimes B) 
\end{equation}
for any two differential graded algebras $A,$ $B.$
Note first that 
because of (\ref{ses:ntb}),
the map $N$ induces an isomorphism 
\begin{equation}  \label{ses:mu}
\overline{C}_{\bullet}(A) \simeq (\text{Ker}(\operatorname{id} - \tau), b')
\end{equation}
where in the right hand side $\operatorname{id}-\tau$ is considered as an operator on $\overline{A}^{\bullet + 1}.$ 

\begin{lemma}
If one identifies $\overline{C}^{\lambda}_{\bullet}$ with the right hand side of (\ref{ses:mu}), then the shuffle product is a morphism of complexes 
\begin{equation} \label{ses:bujumbura}
\times: \overline{C}^{\lambda}_{\bullet}(A) \otimes \overline{C}^{\lambda}_{\bullet}(B) @>>> \overline{C}^{\lambda}_{\bullet + 1}(A \otimes B) 
\end{equation}
\end{lemma}

The same formula defines a product 
\begin{equation} \label{ses:prodred1}
\times: \overline{C}^{\lambda}_{\bullet}(A[\eta]) \otimes \overline{C}^{\lambda}_{\bullet}(B[\eta]) @>>> \overline{C}^{\lambda}_{\bullet + 1}((A \otimes B) [\eta])
\end{equation}

\subsection{Operations on the cyclic complexes} \label{ss:opera}
Here we recall, in a modified version, some results from \cite{NT2}. Let $u$ be a formal parameter of degree two. Consider  the differential graded algebra $A[\eta] = A + A\eta,$ $\text{deg} \eta = -1$ (unlike in Lemma \ref{lemma:ccredeta}), $\eta ^2 = 0$ with the differential $\frac{\partial}{\partial \eta}.$ 
Consider the complex $\overline{C}^{\lambda}_{\bullet}(A[\eta])[u]$ with the differential $\frac{\partial}{\partial \eta} + u\cdot b.$
\begin{thm} \label{thm:calc3}
There exist natural pairings of $k[[u]]$-modules 
$$\bullet: \overline{C}^{\lambda}_{\bullet - 1}(A[\eta])[[u]] \otimes CC^-_{- \bullet}(A) @>>>
CC^-_{- \bullet}(A)
$$
$$\bullet: \overline{C}^{\lambda}_{\bullet - 1}(A[\eta])[[u, u^{-1}] \otimes CC^{\text{per}}_{- \bullet}(A) @>>>
CC^{\text{per}}_{- \bullet}(A)
$$
such that: 
\begin{enumerate}
\item $ \eta ^{(m)} \bullet ? = \text{Id}$ for $m>0.$
\item For $x_i \in A$ the operation $(x_1 \otimes \cdots \otimes x_p) \bullet ?$ sends $C_N (A)$ to \newline
$\sum_{i,j \geq 0} C_{N-p+i}(A)u^j$.
\item The component of  $(x_1 \otimes \cdots \otimes x_p) \bullet (a_0 \otimes \cdots \otimes a_N)$ in $C_{N-p}[[u]]$ is equal to 
$$ \frac{1}{p!} \sum_{i=1}^{p} (-1)^{i(p-1)}a_0 [x_{i+1}, a_1] [x_{i+2}, a_2]\cdots [x_{i}, a_p] \otimes a_{p+1} \otimes \cdots \otimes a_N
$$
\item $(x_1 \otimes \cdots \otimes x_p) \bullet 1 =  \sum_{i=1}^{p} (-1)^{i(p-1)} 1 \otimes x_{i+1} \otimes x_{i+2} \otimes \cdots \otimes x_i.$
\end{enumerate}
\end{thm}

This theorem follows from two others which we include for the readers convenience. Let $(C^{\bullet}(A,A), \delta)$ be the Hochschild cochain complex of $A$. We denote its underlying differential graded algebra with the cup product by ${\cal {E}}^{\bullet}_A$ and its underlying differential graded Lie algebra with the Gerstenhaber bracket by ${{\frak{g}}}^{\bullet - 1}_A.$ Recall that the Hochschild complex $C_{\bullet}(A,A)$ is a module over the differential graded Lie algebra ${{\frak{g}}}^{\bullet}_A ;$ we denote the action of a cochain $D \in {{\frak{g}}}^{\bullet}_A $ on $C_{\bullet}(A,A)$ by $L_D$. At the level of cohomology, if $A$ is the ring of functions on a smooth manifolds, this action induces the action of multivectors on forms via the Lie derivative $L_D = [d, \iota_D].$

Recall that the complex $CC^-_{\bullet}(A)$ (or $CC^{\text{per}}_{\bullet}(A)$) carries a multiplication (\cite{HJ}) which is a morphism of complexes if $A$ is commutative.

\begin{thm} \label{thm:calc1}
There exist natural morphisms of complexes 
$$
\bullet : CC^-_{\bullet}(A) \otimes CC^-_{\bullet}({\cal{E}}^{\bullet}_A) @>>> CC^-_{\bullet}(A)
$$
$$
\bullet : CC^{\text{per}}_{\bullet}(A) \otimes CC^{\text{per}}_{\bullet}({\cal{E}}^{\bullet}_A) @>>> CC^{\text{per}}_{\bullet}(A)
$$

such that:
\begin{enumerate}
\item For a Hochschild cochain $D \in C^p(A,A)$
$$(a_0 \otimes \cdots \otimes a_N) \bullet D = (-1)^{pN}a_0D(a_1, \cdots, a_p) \otimes a_{p+1} \cdots \otimes a_N
$$
$$(a_0 \otimes \cdots \otimes a_N) \bullet (1 \otimes D) = (-1)^{(p-1)N}L_D(a_0 \otimes \cdots \otimes a_N)
$$
\item If one views  $CC^-_{\bullet}(A)$ as a subspace of  $CC^-_{\bullet}({\cal{E}}^{\bullet}_A)$, then the induced product on  $CC^-_{\bullet}(A)$ becomes the Hood-Jones product from \cite{HJ}; the same for $CC^{\text{per}}.$
\end{enumerate}
\end{thm}

The product $\bullet$ is given by the formula  
$$         
     a\bullet x =  a\bullet _1 x +  a\bullet _2 x + u \cdot a\bullet_3 x
$$
where $\bullet _i$ are constructed in \cite{NT2}. (One can show that in fact $CC^-_{\bullet}({\cal{E}}^{\bullet}_A)$ is an $A_{\infty}$ algebra and $ CC^-_{\bullet}(A)$ is an  $A_{\infty}$ module over this algebra (this structure extends the  $A_{\infty}$ structure from \cite{GJ})).

Let ${{\frak{g}}}^{\bullet}_A[\epsilon] = {{\frak{g}}}^{\bullet}_A + {{\frak{g}}}^{\bullet}_A\epsilon,$ $\epsilon ^2 = 0$, $\text{deg} \epsilon = 1.$
Let $u$ be a formal parameter of degree 2 and consider the graded algebra $U({{\frak{g}}}^{\bullet}_A[\epsilon])[u]$ with the differential $u \cdot \frac{\partial}{\partial \epsilon} + \delta$ (recall that $\delta$ is the differential in ${{\frak{g}}}^{\bullet}_A$). Let $\star$ be the Hood - Jones product on  $CC^-_{\bullet}({\cal{E}}^{\bullet}_A).$
\begin{thm} \label{thm:calc2}
There exist morphisms of complexes of $k[u]$-modules
\begin{eqnarray*}
U({{\frak{g}}}^{\bullet}_A[\epsilon])[u] \otimes _{U({{\frak{g}}}^{\bullet}_A)[u]} CC^-_{\bullet}(A) & @>>> & CC^-_{\bullet}(A) \\
U({{\frak{g}}}^{\bullet}_A[\epsilon])[u,u^{-1}] \otimes _{U({{\frak{g}}}^{\bullet}_A)[u,u^{-1}]} CC^{\text{per}}_{\bullet}(A) & @>>> & CC^{\text{per}}_{\bullet}(A)
\end{eqnarray*}
such that
\begin{multline*}
(\epsilon  D_1\cdot \dots \cdot \epsilon D_m) \bullet \alpha = I(D_1, \ldots, D_m)\alpha = \\
(-1)^{\deg \alpha\sum_{i=1}^{m}\text{deg}D_i}\frac{1}{m!}\underset{\sigma}\sum
\epsilon_{\sigma} \alpha \bullet (D_{\sigma _1} \star (D_{\sigma _2} \star (\ldots \star D_{\sigma _m}))\ldots )
\end{multline*}
\end{thm}
The motivation for this theorem is the following. If $A=C^{\infty}(M)$ then $HC^-_{\bullet}(A)$ is isomorphic to the cohomology of the complex $\Omega^{\bullet}(M)[[u]], ud;$ the cohomology of ${{\frak{g}}}^{\bullet}_A$ is the graded Lie algebra ${{\frak{g}}}^{\bullet}_M$ of multivector fields on $M$; one can define an action of ${{\frak{g}}}^{\bullet}_M$ on $HC^-_{\bullet}(A)$: for two multivector fields $X,$ $Y$ the action of $X+\epsilon Y$ is given by $L_X+\iota_Y$ where $\iota_Y$ is the contraction operator and $L_X=[d,\iota_X].$

{\bf {Remark.}} One can extend the pairings in theorem \ref{thm:calc2} to an action up to homotopy of ${{\frak{g}}}^{\bullet}_A[\epsilon]$ on $CC^-_{\bullet}(A)$, i.e., to an $L_{\infty}$ morphism from ${{\frak{g}}}^{\bullet}_A[\epsilon]$ to $\text{End}^{\bullet}_k ( CC^-_{\bullet}(A)).$

For the readers convenience let us show how to reduce Theorem \ref{thm:calc3} from Theorem \ref{thm:calc2}. Let $M(A)$ be the algebra of matrices
$(a_{ij})_{0\leq i, j\le\infty}$ for which $ a_{ij}\in A$ and all but finitely many of
$a_{ij}$ are zero. The same space considered as a Lie algebra is
${{\frak{gl}}}(A)$.

 Let ${M}_\infty(A)$ be the associative algebra of matrices
$(a_{ij})_{1\leq i, j\le\infty}$ such that $a_{ij}\in A$ and $a_{ij}=0$
on all but finitely many diagonals. There is a morphism of differential graded Lie algebras
${{\frak{gl}}}(A[\eta]) @>>> {{\frak{g}}}^{\bullet}_A$ where $\eta$ is a formal
parameter of degree $-1$ such that $\eta^2=0$; the differential
on $A[\eta]$ is $\partial/\partial\eta$. This morphism sends $\eta a$
to $a$ in $C^0$ and $a$ to $\text {ad}(a)$ in $C^1$.

Consider the following morphisms:
$$
CC_{\bullet}^{per}(A) @>i>> CC_{\bullet}^{per}({M}_\infty(A));
$$
$$
CC_{\bullet}^{per}({M}(A)) @>\text{tr}>> CC_{\bullet}^{per}(A)  
$$
The map $i$ is induced by the inclusion $a\mapsto a\cdot 1$; the second
map $tr$ acts as follows:
$$
(a_0 m_0,\dots, a_n m_n) \mapsto
tr(m_0\dots m_n)\; (a_0,\dots, a_n)
$$
for $a_i\in A$ and $m_i\in {M}(k)$. One checks easily (\cite{NT2}) that the restriction of the pairing from Theorem \ref{thm:calc2} to ${{\frak{gl}}}(A[\eta])$, when composed with $i$ and $tr$, descends to a pairing
$$
(k[u] \otimes _{U({{\frak{gl}}} (k[\eta]))[u]} U({{\frak{gl}}}(A[\epsilon, \eta]))[u] \otimes _ {U({{\frak{gl}}}(k[\epsilon]))[u]}) \otimes CC^-_{- \bullet}(A) @>>>  CC^-_{- \bullet}(A)
$$  
Finally, consider the complex $C_{\bullet}({{\frak{gl}}}(A[\eta]), {{\frak{gl}}}(k))[u]$ with the differential $\partial + u\cdot \frac{\partial}{\partial \eta}.$  Let ${{\frak{g}}}={{\frak{gl}}}(A)$ and ${{\frak{h}}}={{\frak{gl}}}(k).$
One constructs a morphism of complexes 
$$
C_{\bullet}({{\frak{g}}}[\eta], {{\frak{h}}})[u] @>>> k[u] \otimes _{U({{\frak{h}}} [\eta])[u]} U({{\frak{g}}}[\epsilon, \eta])[u] \otimes _ {U({{\frak{h}}}[\epsilon])[u]} k[u]
$$
as follows. To describe the image of a chain
$D_1\wedge\dots\wedge D_n \wedge \eta E_1\wedge\dots\wedge\eta E_m$,
write the expression
$$
D_1(\epsilon-\eta)\cdot\dots\cdot D_n(\epsilon-\eta)\cdot E_1\epsilon\eta
\cdot\dots\cdot E_m\epsilon\eta
$$
and then represent it as a sum
$$
\sum\pm D_{\iota_1}\eta\cdot\dots\cdot D_{\iota_k}\eta\cdot
E_1\epsilon\eta\cdot\dots\cdot E_n\epsilon\eta\cdot
D_{j_1}\epsilon\cdot\dots\cdot D_{j_m}\epsilon
$$
in the symmetric algebra of the graded space ${\frak{g}}[\epsilon,\eta]$.
For example $D\mapsto D\epsilon-D\eta$;
$$
D_1\wedge D_2\longmapsto
D_1\eta\cdot D_2\eta+D_1\eta\cdot D_2\epsilon+
(-1)^{(|D_1|+1)(|D_2|+1)}D_2\eta\cdot D_1\epsilon+
D_1\epsilon\cdot D_2\epsilon 
$$
$$
D_1\wedge D_2\eta\longmapsto-D_1\eta\cdot D_2\epsilon\eta+
(-1)^{(|D_1|+1)(|D_2|+1)}D_2\epsilon\eta\cdot D_1\eta
$$
etc.

Thus one gets a pairing of $C_{\bullet - 1}({{\frak{gl}}}(A[\eta]), {{\frak{gl}}}(k))[u]$ with $CC^-_{- \bullet}(A);$ Theorem \ref{thm:calc3} now follows from Theorem \ref{thm:liecyclic}.
\subsection{Rigidity of periodic cyclic homology} \label{ss:rigid}
Let $A_0$ be an associative algebra over a ring $k$ of characteristic zero. Let $t$ be a formal parameter. Consider a formal deformation of $A_0,$ that is, an associative $k[[t]]$-linear product $*$ on $A[[t]]$ such that $a*b = ab + O(t)$ for $a,b \in A_0$ and $1*a=a*1=a.$ We will denote by $A$ the algebra $(A_0[[t]],*)$ and by $A_0[[t]]$ the (undeformed) algebra of power series with values in $A_0.$

There are two versions of Hochschild, cyclic, etc. complexes for complete $k[[t]]$-algebras: one is defined using tensor products over $k$, the other using completed tensor products over $k[[t]].$ We will denote corresponding periodic cyclic complexes by $CC_{\bullet}^{\text{per}}(-)_k ,$resp. by complexes by  $CC_{\bullet}^{\text{per}}(-)_{k[[t]]}.$ Note that there are obvious projection maps 
$p: CC_{\bullet}^{\text{per}}(-)_{k} @>>> CC_{\bullet}^{\text{per}}(-)_{k[[t]]}.$ Define also the map $\sigma : CC_{\bullet}^{\text{per}}(A)_{k} @>>> CC_{\bullet}^{\text{per}}(A_0).$
\begin{thm} \label{thm:rigid} (Goodwillie)
\begin{enumerate}
\item The map $\sigma$ is a quasi-isomorphism.
\item There is a canonical isomorphism of complexes $\chi$ making the following diagram commute:
$$
\begin{array}{ccc}
CC_{\bullet}^{\text{per}}(A)_k & @>\sigma>> & CC_{\bullet}^{\text{per}}(A_0) \\
\downarrow{p} & & \downarrow{p} \\
CC_{\bullet}^{\text{per}}(A)_{k[[t]]} & @>\chi>> & CC_{\bullet}^{\text{per}}(A_0[[t]])_{k[[t]]}
\end{array}
$$
\end{enumerate}
\end{thm}
\begin{pf}
First, let us construct $\chi.$ Put $\lambda(a_1,a_2)=a_1a_2 - a_1 *a_2;$ then $\lambda$ is an element of ${{\frak{g}}}^1_A$ such that 
$$\delta\lambda + \frac{1}{2}[\lambda, \lambda] = 0.$$
Since $\lambda$ is divisible by $t,$ the operation 
\begin{equation} \label{ses:chi}
\chi = \sum_{m=0}^{\infty}\frac{1}{m!}I(\lambda, \ldots, \lambda)
\end{equation}
is well-defined on $CC_{\bullet}^{\text{per}}(A)_{k[[t]]}$ (where $I(\lambda, \ldots, \lambda)$ are the operators from theorem \ref{thm:calc2}). One checks that $\chi$ is indeed a morphism of complexes (\cite{NT2}).It is equal to identity modulo $t,$ therefore it is invertible.

Note that formula (\ref{ses:chi}) also defines an isomorphism
$$\chi _k : CC_{\bullet}^{\text{per}}(A)_k @>>> CC_{\bullet}^{\text{per}}(A_0[[t]])_k
$$
Now let us observe that the maps 
$$
CC_{\bullet}^{\text{per}}(A_0[[t]])_k @>\sigma>> CC_{\bullet}^{\text{per}}(A_0)
$$
and
$$
CC_{\bullet}^{\text{per}}(A_0) @>i>> CC_{\bullet}^{\text{per}}(A_0[[t]])
$$
induce mutually inverse isomorphisms on homology. Indeed, the operator $L_{\frac{\partial}{\partial t}}$ acts by zero on homology because 
$$
L_{\frac{\partial}{\partial t}} = [b+uB, I({\frac{\partial}{\partial t}})]
$$
Thus we get a commutative diagram in which all horizontal maps are quasi-isomorphisms:
$$
\begin{array}{ccccc}
CC_{\bullet}^{\text{per}}(A)_k & @>\chi _k >> & CC_{\bullet}^{\text{per}}(A_0[[t]])_k& @>\sigma>> & CC_{\bullet}^{\text{per}}(A_0) \\
\downarrow{p} & & \downarrow{p} & & \downarrow{i} \\
CC_{\bullet}^{\text{per}}(A)_{k[[t]]} & @>\chi>> & CC_{\bullet}^{\text{per}}(A_0[[t]])_{k[[t]]} & = & CC_{\bullet}^{\text{per}}(A_0[[t]])_{k[[t]]} 
\end{array}
$$
But $\sigma \chi_k = \sigma,$ which concludes the proof.
\end{pf}
%
%
%
%
%
%
%
%
%
%
%
%

\section{Appendix 2. Review of Lie algebra cohomology and characteristic classes}  \label{section:lie}
For a differential graded Lie algebra (DGLA)
$({\frak{g}}, \delta)$ over $k$ let $U({\frak{g}})$ denote the universal
enveloping algebra of $({\frak{g}}, \delta)$ (this is a DGA with the
differential $\delta$ acting by the Leibnitz rule).

For a DGLA $({\frak{g}}, \delta)$ over $k$ let ${\frak{g}}[\epsilon]
\overset{def}{=}
({\frak{g}},\delta)\otimes_k 
(k[\epsilon]/\epsilon^2, \displaystyle\frac{\partial}{\partial\epsilon})$,
where $\deg\epsilon = 1$.

For a DGLA $({\frak{g}}, \delta)$ over $k$ and a DGL-subalgebra ${\frak{h}}$
let
\begin{eqnarray*}
C_{\bullet}({\frak{g}}) & \overset{def}{=} &
U({\frak{g}}[\epsilon])\otimes_{U({\frak{g}})}k\ , \\
C_{\bullet}({\frak{g}})_{{\frak{h}}} & \overset{def}{=} &
k\otimes_{U({\frak{h}})} U({\frak{g}}[\epsilon])\otimes_{U({\frak{g}})}k\ , \\
C_{\bullet}({\frak{g}},{\frak{h}}) & \overset{def}{=} &
k\otimes_{U({\frak{h}}[\epsilon])} U({\frak{g}}[\epsilon])\otimes_{U({\frak{g}})}k\ .
\end{eqnarray*}
The structure of a DG coalgebra on $U({\frak{g}})$ induces similar structures
on all of the above complexes.

For a left ${\frak{g}}$-module $L$ put

\begin{equation} \label{eq:cochcomplie}
C^{\bullet}({\frak{g}}, L) = \operatorname{Hom}_{U({\frak{g}})}(U({\frak{g}}[\epsilon]), L)
\end{equation}
\begin{equation} \label{eq:cochcomplierel}
C^{\bullet}({\frak{g}}, {\frak{h}}; L) = \operatorname{Hom}_{U({\frak{g}})}(U({\frak{g}}[\epsilon]) \otimes _{U({\frak{h}}[\epsilon])}, L)
\end{equation}
\begin{equation} \label{eq:clieconst}
C^{\bullet}({\frak{g}}) = C^{\bullet}({\frak{g}}, k); \;\; 
C^{\bullet}({\frak{g}}, {\frak{h}}) = C^{\bullet}({\frak{g}}, {\frak{h}}; k)
\end{equation}

The differential in $C^{\bullet}$, $C_{\bullet}$ will be denoted by $\partial$. The complexes  $C^{\bullet}({\frak{g}})$ and $C^{\bullet}({\frak{g}}, {\frak{h}})$ are differential graded commutative algebras; $C^{\bullet}({\frak{g}}, L)$ is a DG module over $C^{\bullet}({\frak{g}}, )$ and $C^{\bullet}({\frak{g}}, {\frak{h}}; L)$ is a DG module over $C^{\bullet}({\frak{g}}, {\frak{h}})$. Of course our definitions are equivalent to the standard ones; the differential $\partial $ consists of two parts, the Chevalley-Eilenberg differential and the one induced by the differential on ${\frak{g}}$.

For an element $X$ of ${\frak{g}}$, right multiplication by $X$ (respectively by $\epsilon X$) on $U({\frak{g}})[\epsilon]$ induces a derivation $L_X$, resp. $\iota _X$ of the DG algebra $C^{\bullet}({\frak{g}})$ and of the DG module $C^{\bullet}({\frak{g}}, L)$. These derivation satisfy
\begin{equation} \label{eq:cartann0}
[\partial , \iota_X]=L_X;\:\:\:
[L_X,L_Y]=L_{[X,Y]};\:\:\:[\partial ,L_X]=0;
\end{equation}
$$ \:\:\: [L_X,\iota_Y]=\iota_{[X,Y]};\:\:\:[\iota_X,\iota_Y]=0$$
for $X,Y$ in $\g$.

%
%
%
%
%
%
%
%
%
%
%
%
%
%
\subsection{Characteristic classes with coefficients in the trivial module}
\label{chclasstrivco}
In this subsection we will construct the map $S^m({\frak{h}} ')^{{\frak{h}}}\rightarrow
H^{2m}({\frak{g}},{\frak{h}})$ from the space of invariant polynomials on ${\frak{h}}$ to the relative
Lie algebra cohomology. We will first construct it at the level of cohomology
and then represent by explicit cochains.

Let $\eta$ be a formal parameter of degree $1$, $\eta ^2=0$. Construct the
differential graded Lie algebra $(\g[\eta],\frac{\partial}{\partial \eta})$ exactly as in the beginning of Section \ref{section:lie}.
Consider the map 
\begin{equation}   \label{eq:gtogeta}
C_{\bullet}(\g,{\frak{h}})\rightarrow C_{\bullet}(\g[\eta],{\frak{h}}) 
\end{equation}
\begin{lemma}    \label{getaheta}
i) The inclusion
$$C_{\bullet}({\frak{h}}[\eta],{\frak{h}}) \rightarrow C_{\bullet}(\g[\eta],{\frak{h}})$$
is a quasi-isomorphism.

ii) One has
$$C_{\bullet}({\frak{h}}[\eta],{\frak{h}})=H_{\bullet}({\frak{h}}[\eta],{\frak{h}})=S^{\frac{1}{2}{\bullet}}({\frak{h}})_{{\frak{h}}}.$$
\end{lemma}
We leave the proof to the reader.

We get a map $H_{\bullet}(\g,{\frak{h}})\stackrel{\alpha}{\rightarrow} S^{\frac{1}{2}\bullet}({\frak{h}})_{{\frak{h}}}$; for any
invariant polynomial $P$ from $S^m({\frak{h}} ')^{{\frak{h}}}$ put
\begin{equation} \label{eq:chern0}
c_P=m!P \circ \alpha \in H^{2m}(\g, {\frak{h}}).
\end{equation}
Our next aim is 
to represent these cohomology classes by explicit cochains.

 Choose an ${\frak{h}}$-module decomposition $\g={\frak{h}}\oplus V$; define the cochain $A$
in $C^1(\g)\otimes {\frak{h}}$ where ${\frak{h}}$ is viewed as trivial $\g$-module as follows:
\begin{equation} \label{eq:relialconn}
A(X)=(Projection\: of\: X\: to\: {\frak{h}} \: along\: V)
\end{equation}
For ${\frak{h}}$-valued forms on $\g$ (or, more generally, for elements of $C^{\bullet}\otimes
{\frak{h}}$ where $C^{\bullet}$ is a commutative differential graded algebra) we define the bracket 
\begin{equation}  \label{eq:brach}
[\varphi _1 \otimes h_1, \varphi _2 \otimes h_2]=(\varphi _1 \cdot \varphi _2)
\otimes [h_1,h_2]
\end{equation}
for $\varphi _i$ in $C^{\bullet}$ and $h_i$ in ${\frak{h}}$. The complex $(C^{\bullet}(\g)
\otimes {\frak{h}} , \partial ^{Lie})$ becomes a differential graded Lie algebra.
One has
\begin{equation}  \label{eq:propa}
L_hA+[h,A]=0;\:\:\: \iota_hA=h
\end{equation}
for $h$ in ${\frak{h}}$.
 Put
\begin{equation}  \label{eq:curvature}
R=\partial ^{Lie} A+\frac{1}{2} [A,A]
\end{equation}
One has 
\begin{equation}  \label{eq:curvature1}
R(X,Y)=[A(X),A(Y)]-A([X,Y])
\end{equation}
for $X,Y \in \g$. 
\begin{lemma}  \label{propr}
For any $h$ in ${\frak{h}}$
$$\iota_hR=0;\:\:\:L_hR+[h,R]=0;\:\:\:\partial ^{Lie}R+[A,R]=0$$
\end{lemma}
{\bf Proof.} Follows from (\ref{eq:propa}) and \ref{eq:curvature} by a
straightforward computation. $\Box$

Note that any multilinear form on ${\frak{h}}$ extends to a $C^{\bullet}(\g)$-valued form on
$C^{\bullet}(\g) \otimes {\frak{h}}$: just put
$$P(\varphi _1 \otimes h_1, \ldots , \varphi _m\otimes h_m)=\varphi _1 \ldots
\varphi _m \otimes P(h_1 , \ldots , h_m)$$
\begin{definition}  \label{chern}
{\em Let $P$ be an element of $S^m({\frak{h}} ')^{{\frak{h}}}$. Put
$$c_P=P(R,\ldots, R) \in C^{2m}(\g)$$}
\end{definition}
\begin{proposition} \label{procher}
i) The cochain $c_P$ is a cocycle in $C^{2m}(\g , {\frak{h}}).$

ii) The cohomology class of this cocycle is given by the formula
(\ref{eq:chern0}). In particular it does not depend on a choice of a decomposition
$\g={\frak{h}}\oplus V$.
\end{proposition}
{\bf Proof.} Using Lemma \ref{propr} one sees that
$$\partial ^{Lie}c_P = \sum_{i=1}^{m} P(R, \ldots , [A,R], \ldots, R)=0$$
because $P$ is ${\frak{h}}$-invariant; by the same Lemma, $\iota_hc_P=0$ and 
$$L_h c_P = \sum_{i=1}^{m} P(R, \ldots , [h,R], \ldots, R)=0.$$ This proves
i). To prove ii), note that the cochain $A$ extends to the differential graded
Lie algebra $\g[\eta]$ by putting 
$$A(X+\eta Y)=A(X).$$ Put 
\begin{equation}  \label{eq:dgcurvature}
R=\partial ^{Lie} A+\frac{\partial}{\partial \eta}A +\frac{1}{2} [A,A].
\end{equation}
Then all the properties of $A$ and $R$ (formulas \ref{eq:propa}, Lemma
\ref{propr}, the statement i)) still hold if one replaces $\partial
^{Lie}$ by $\partial ^{Lie} +\frac{\partial}{\partial \eta}$. Now note that 
\begin{equation}  \label{eq:ceta}
R(\eta h)=h;\:\:\:P(R, \ldots , R)(\eta h, \ldots , \eta h)=m!P(h, \ldots , h)
\end{equation}
which is the same as in (\ref{eq:chern0}).  $\Box$
\subsection{The Weil algebra}   \label{theweilalg}
\label{weilalg}
\begin{em} 
The {\em Weil algebra} of a Lie algebra ${\frak{h}}$ is the relative cochain complex
$C^{\bullet}({\frak{h}}[\eta],{\frak{h}})$ with the exterior product and
the differential $$\partial^{W}=\partial^{Lie}+\frac{\partial}{\partial
\eta}$$
\end{em}
Denote this differential graded Lie algebra by $W^{\bullet}({\frak{h}})$.
\begin{remark}   \label{weilalg1}
\begin{em}
Let us give an equivalent definition. Let $W^{\bullet}({\frak{h}})$ be the free commutative
differential graded algebra with unit generated by 
$$W^1({\frak{h}})={\frak{h}} ',$$
$$W^2({\frak{h}})={\frak{h}} '$$
(${\frak{h}}'$ is the dual space to ${\frak{h}}$); by $A,$ resp. $R$ we denote the element of
$W^1({\frak{h}})\otimes {\frak{h}},$ resp. $W^2({\frak{h}})\otimes {\frak{h}},$ corresponding to $id_{{\frak{h}}};$
then the differential $\partial^W$ is by definition given by 
\begin{equation} \label{eq:wealdi}
\partial^W A=R-\frac{1}{2}[A,A];\:\:\:\partial^W R = - [A,R]
\end{equation}
The bracket in the above formula is given by (\ref{eq:brach})
\end{em}
\end{remark}
Let $h$ be an element of ${\frak{h}}$. For $p > 0$ define $\iota_h: W^{p}({\frak{h}}) \rightarrow W^{p-1}({\frak{h}})$
to be the unique derivation which, for $l \in {\frak{h}}' ,$ sends the corresponding
element of $W^1$ to $l(h) \cdot 1$ and the corresponding element of $W^2$ to
zero. Let $L_h :W^{\bullet}({\frak{h}})\rightarrow W^{\bullet}({\frak{h}})$ be the unique derivation whose
restriction to $W^1$ and $W^2$ corresponds to the operator $Ad^*_h$ on ${\frak{h}}'.$
The operators $L_h$ and $\iota_h$ are the same as those defined before (\ref{eq:cartann0}). Either using this or checking directly one gets
\begin{equation}  \label{eq:cartan}
[\partial^{W}, \iota_X]=L_X;\:\:\: [L_X,L_Y]=L_{[X,Y]};\:\:\:[\partial^{W},L_X]=0;
\end{equation}
$$\:\:\: [L_X,\iota_Y]=\iota_{[X,Y]};\:\:\:[\iota_X,\iota_Y]=0$$
for $X$, $Y$ from ${\frak{h}}$.
One has also 
\begin{equation}   \label{eq:ia,la}
\iota_h A = h; \:\:\: L_h A + [h,A] = 0
\end{equation}
\begin{definition}   \label{bas}
\begin{em}
Let $(W^{\bullet}, \partial^{W})$ be any complex equipped with operators $\iota_X$ of
degree $-1$ and $L_X$ of degree zero which are linear in $X \in {\frak{h}}$ and
satisfy (\ref{eq:cartan}). We call an element $\varphi$ of $W^{\bullet}$ {\em basic} if $L_h \varphi = \iota_h
\varphi = 0$ for any $h \in {\frak{h}}$. The subcomplex of basic elements of $W^{\bullet}$ is
denoted by $W^{\bullet}_{bas}.$
\end{em}
\end{definition}
We can now give another definition of the characteristic classes $c_P$. Let
$\g$ be a Lie algebra, ${\frak{h}}$ its subalgebra such that $\g$ is completely
reducible with respect to  ${\frak{h}}$. Consider the ${\frak{h}}$ - valued cochains $A$ and
$R$; they define a morphism of differential graded algebras 
\begin{equation}  \label{eq:cwoloch}
W^{\bullet}({\frak{h}}) \rightarrow C^{\bullet}(\g);
\end{equation}
because of (\ref{eq:propa}) and (\ref{eq:ia,la}) this map commutes with the
operators $\iota_h$ and $L_h$ and therefore restricts to the subcomplexes of basic
elements. But 
$$C^{\bullet}(\g)_{bas}=C^{\bullet}(\g,{\frak{h}});\:\:\:\:W^{\bullet}({\frak{h}})_{bas}=S^{\frac{1}{2}*}({\frak{h}}')^{{\frak{h}}};$$
therefore one gets a morphism of differential graded algebras 
\begin{equation}   \label{eq:cw1}
S^{\frac{1}{2}*}({\frak{h}}')^{{\frak{h}}} \rightarrow C^{\bullet}(\g,{\frak{h}})
\end{equation}
Clearly, this map sends an invariant polynomial $P$ on ${\frak{h}}$ to the cocycle 
$c_P$ from definition \ref{procher}.
\subsection{Characteristic classes in relative Lie algebra cohomology with
coefficients in homotopically constant complexes of modules}  \label{ss:charactlie}
Let, as above, $\g$ be a Lie algebra and ${\frak{h}}$ a Lie subalgebra of $\g$ which is reductive
and whose action on $\g$ is completely reducible. Assume that $(L^{\bullet},
\partial)$ is a complex of $\g$-modules on which ${\frak{h}}$ acts completely
reducibly. We denote the action of an element $X$ of $\g$ on $L^{\bullet}$ by $L_X$. 
\begin{definition} \label{hotri}
\begin{em}The complex $L^{\bullet}$ is said to be {\em homotopically
constant} if there exist operators $\iota_X : L^{p} \rightarrow L^{p-1},$ $X \in
\g,$ $p>0$, which are linear in $X$ and satisfy
\begin{equation} \label{eq:cartann}
[\partial , \iota_X]=L_X;\:\:\:
[L_X,L_Y]=L_{[X,Y]};\:\:\:[\partial ,L_X]=0;
\end{equation}
$$ \:\:\: [L_X,\iota_Y]=\iota_{[X,Y]};\:\:\:[\iota_X,\iota_Y]=0$$
for $X,Y$ in $\g$.
\end{em}
\end{definition}
Let $L^{\bullet}$ be homotopically trivial. We shall construct morphisms of complexes 
\begin{equation}  \label{eq:cwolota}
c:W^{\bullet}({\frak{h}})\otimes L^{\bullet} \rightarrow C^{\bullet}(\g, L^{\bullet})
\end{equation}
and
\begin{equation}  \label{eq:cwb}
c:(W^{\bullet}({\frak{h}})\otimes L^{\bullet})_{bas} \rightarrow C^{\bullet}(\g,{\frak{h}}; L^{\bullet}).
\end{equation}
(The operators $\iota_h, \:\: L_h$ on
the left hand side of (\ref{eq:cwolota}) are defined by
$$\iota_h(w\otimes l)=\iota_hw \otimes l + (-1)^{\mid w \mid} w\otimes \iota_h l;$$
$$L_h(w\otimes l)=L_hw \otimes l + w\otimes L_h l).$$
The latter will be the restriction of the former to the subspace of basic
elements. To define the morphisms $c,$ construct the map
\begin{equation}  \label{eq:cwo}
 L^{\bullet} \rightarrow C^{\bullet}(\g, L)
\end{equation}
as follows. For $l \in L^m$ consider the cochains $\lambda ^p \in C^p(\g,
L^{m-p})$;
$$\lambda ^p(X_1, \ldots , X_p)=\iota_{X_p} \ldots \iota_{X_1}(l)$$
Then the map (\ref{eq:cwo}) sends an element $l$ of $L^{\bullet}$ to the cochain
\begin{equation}  \label{eq:cwoo}
\varphi_{l} =  \sum_{p \geq 0}{\lambda ^p}
\end{equation}
Consider also the morphism (\ref{eq:cwolota}). We get the map
\begin{equation} \label{eq:cwol}
W^{\bullet}({\frak{h}})\otimes L^{\bullet} \rightarrow C^{\bullet}(\g)\otimes L^{\bullet} \rightarrow C^{\bullet}(\g, L^{\bullet})
\end{equation}
(the last arrow denotes the $\smile$ product).
\begin{lemma}   \label{lcw}
i) For $h \in {\frak{h}}$ and $l \in L^{\bullet}$ one has
$$\varphi _{\iota_h l}=\iota_h \varphi _l$$
where in the right hand side $\iota_h$ and $L_h$ are those defined before (\ref{eq:cartann0}).

ii) For $l \in L^{\bullet}$ and $a \in U(\g[\varepsilon])$ one has
$$\varphi _{\partial l}(a)= \partial ^{Lie} \varphi_l(a) + (-1)^{\mid a
\mid}\partial \varphi_l (a)$$
In other words, the map (\ref{eq:cwo}) is a morphism of complexes and it
commutes with the operators $L_h$, $\iota_h$ for $h \in {\frak{h}}$.
\end{lemma}
{\bf Proof} The first equation of i) is obvious. The second follows from the
first, from ii) and from $[\partial, \iota_h]=L_h$,  $[\partial^{Lie},
\iota_h]=L_h$. To prove ii), note that formulas (\ref{eq:cartann}) mean that the
differential graded algebra $(U(\g [\varepsilon]),
 \frac{\partial}{\partial \varepsilon})$ acts on $(L^{\bullet}, \partial)$. For $a \in U(\g [\varepsilon]$ put
$\tilde{a}= (-1)^{\mid a \mid (\mid a \mid - 1)/2}a;$ then
$$(\partial^{Lie} \varphi _l)(a)=\varphi _l( \frac{\partial a}{\partial
\varepsilon})=( \frac{\partial a}{\partial \varepsilon}){\tilde{} } \cdot
l=[\tilde{a}, \partial] \cdot l = $$
$$\tilde{a}\partial l - (-1)^{\mid a
\mid}\partial(\tilde{a}\cdot l)=\varphi_{\partial l}(a)-(-1)^{\mid a
\mid}\partial \varphi _l (a) \:\:\:\: \Box$$
\begin{proposition} \label{cwcwcw}
The composition (\ref{eq:cwol}) is a morphism of complexes commuting with
$\iota_h,\:\:\:L_h$ for $h \in {\frak{h}}$.
\end{proposition}
{\bf Proof} Follows from Lemma \ref{lcw} $\Box$.

\end{document}